\documentclass[a4paper,12pt,titlepage]{article}
\usepackage{latexsym}
\usepackage{amsbsy}
\usepackage{amssymb}
\usepackage{amscd}
\usepackage{amsmath}
\begin{document}

\newcommand{\ep}{\hspace*{\fill}$\Box$}
\newcommand{\eps}{\varepsilon}
\newcommand{\pr}{{\bf Proof. }}
\newcommand{\ms}{\medskip\\}
\newcommand{\cl}{\mbox{\rm cl}}
\newcommand{\g}{\ensuremath{\mathfrak g} }
\newcommand{\gc}{${\cal G}$-complete }
\newcommand{\sa}{\stackrel{\scriptstyle s}{\approx}}
\newcommand{\prol}{\mbox{\rm pr}^{(n)}}
\newcommand{\prolo}{\mbox{\rm pr}^{(1)}}
\newcommand{\deta}{\frac{d}{d \eta}{\Big\vert}_{_{0}}}
\newcommand{\detas}{\frac{d}{d \eta}{\big\vert}_{_{0}}}
\newcommand{\R}{\mathbb R}
\newcommand{\N}{\mathbb N}
\newcommand{\C}{\mathbb C}
\newcommand{\Z}{\mathbb Z}
\newcommand{\K}{\mathbb K}
\newcommand{\sR}{\mathbb R}
\newcommand{\sN}{\mathbb N}
\newcommand{\gK}{{\cal K}}
\newcommand{\gR}{{\cal R}}
\newcommand{\gC}{{\cal C}}
\newcommand{\Dp}{${\cal D}'$ }                                          
\newcommand{\go}{${\cal G}(\Omega)$ }
\newcommand{\grn}{${\cal G}(\R^n)$ }
\newcommand{\grp}{${\cal G}(\R^p)$ }
\newcommand{\grq}{${\cal G}(\R^q)$ }
\newcommand{\gt}{${\cal G}_\tau$ }
\newcommand{\gto}{${\cal G}_\tau(\Omega)$ }
\newcommand{\gtrn}{${\cal G}_\tau(\R^n)$ }
\newcommand{\gtrp}{${\cal G}_\tau(\R^p)$ }
\newcommand{\gtrq}{${\cal G}_\tau(\R^q)$ }
\newcommand{\gtn}{{\cal G}_\tau(\R^n) }
\newtheorem{thr}{\hspace*{-1.1mm}}[section]
\newcommand{\bt}{\begin{thr} {\bf Theorem. }}
\newcommand{\et}{\end{thr}}
\newcommand{\bp}{\begin{thr} {\bf Proposition. }}
\newcommand{\bc}{\begin{thr} {\bf Corollary. }}
\newcommand{\blem}{\begin{thr} {\bf Lemma. }}
\newcommand{\bex}{\begin{thr} {\bf Example. }\rm}
\newcommand{\bexs}{\begin{thr} {\bf Examples. }\rm}
\newcommand{\bd}{\begin{thr} {\bf Definition. }}
\newcommand{\beast}{\begin{eqnarray*}}
\newcommand{\eeast}{\end{eqnarray*}}
\newcommand{\wsc}[1]{\overline{#1}^{wsc}}
\newcommand{\todo}[1]{$\clubsuit$\ {\tt #1}\ $\clubsuit$}
\newcommand{\rem}[1]{\vadjust{\rlap{\kern\hsize\thinspace\vbox%
                       to0pt{\hbox{${}_\clubsuit${\small\tt #1}}\vss}}}}
\newcommand{\ahat}{\ensuremath{\hat{\mathcal{A}}_0(M)} }
\newcommand{\atil}{\ensuremath{\tilde{\mathcal{A}}_0(M)} }
\newcommand{\aqtil}{\ensuremath{\tilde{\mathcal{A}}_q(M)} } 
\newcommand{\ehat}{\ensuremath{\hat{\mathcal{E}}(M)} } 
\newcommand{\emhat}{\ensuremath{\hat{\mathcal{E}}_m(M)} }
\newcommand{\nhat}{\ensuremath{\hat{\mathcal{N}}(M)} }
\newcommand{\ghat}{\ensuremath{\hat{\mathcal{G}}(M)} } 
\newcommand{\lhat}{\ensuremath{\hat{L}_X} }                    
\newcommand{\comp}{\subset\subset}
\newcommand{\al}{\alpha}
\newcommand{\bet}{\beta} 
\newcommand{\ga}{\gamma}
\newcommand{\Om}{\Omega}\newcommand{\Ga}{\Gamma}\newcommand{\om}{\omega}
\newcommand{\si}{\sigma}\newcommand{\la}{\lambda}
\newcommand{\de}{\delta}
\newcommand{\vphi}{\varphi}\newcommand{\dl}{{\displaystyle \lim_{\eta>0}}\,}
\newcommand{\intl}{\int\limits}\newcommand{\su}{\sum\limits_{i=1}^2}
\newcommand{\D}{{\cal D}}\newcommand{\Vol}{\mbox{Vol\,}}
\newcommand{\Or}{\mbox{Or}}\newcommand{\sign}{\mbox{sign}}
\newcommand{\na}{\nabla}\newcommand{\pa}{\partial}
\newcommand{\ti}{\tilde}\newcommand{\T}{{\cal T}} \newcommand{\G}{{\cal G}}
\newcommand{\DD}{{\cal D}}\newcommand{\X}{{\cal X}}\newcommand{\E}{{\cal E}} 
\newcommand{\CC}{{\cal C}}\newcommand{\vo}{\Vol}
\newcommand{\bat}{\bar t}
\newcommand{\bx}{\bar x}
\newcommand{\by}{\bar y} \newcommand{\bz}{\bar z}\newcommand{\br}{\bar r}
\newcommand{\fr}{\frac{1}}\newcommand{\il}{\int\limits}
\newcommand{\nn}{\nonumber}
\newcommand{\supp}{\mathop{\mathrm{supp}}}

\newcommand{\vp}{\mbox{vp}\frac{1}{x}}\newcommand{\A}{{\cal A}}
\newcommand{\Ll}{L_{\mbox{\small loc}}}\newcommand{\Hl}{H_{\mbox{\small loc}}}
\newcommand{\Lll}{L_{\mbox{\scriptsize loc}}}
\newcommand{\be}{ \begin{equation} }\newcommand{\ee}{\end{equation} }
\newcommand{\beq}{ \begin{equation} }\newcommand{\eeq}{\end{equation} }
\newcommand{\bea}{\begin{eqnarray}}\newcommand{\eea}{\end{eqnarray}}
\newcommand{\beas}{\begin{eqnarray*}}\newcommand{\eeas}{\end{eqnarray*}}
\newcommand{\beqs}{\begin{equation*}}\newcommand{\eeqs}{\end{equation*}}
\newcommand{\lb}{\label}\newcommand{\rf}{\ref}
\newcommand{\GL}{\mbox{GL}}\newcommand{\bfs}{\boldsymbol}
\newcommand{\ben}{\begin{enumerate}}\newcommand{\een}{\end{enumerate}}
\newcommand{\ba}{\begin{array}}\newcommand{\ea}{\end{array}}
\newtheorem{thi}{\hspace*{-1.1mm}}[section]
\newcommand{\bthm}{\begin{thr} {\bf Theorem. }}
\newcommand{\bprop}{\begin{thr} {\bf Proposition. }}
\newcommand{\bcor}{\begin{thr} {\bf Corollary. }}
\newcommand{\bdef}{\begin{thr} {\bf Definition. }}
\newcommand{\brem}{\begin{thr} {\bf Remark. }\rm}
\newcommand{\bth}{\begin{thr}\rm}
\newcommand{\ethi}{\end{thr}}
\newcommand{\ca}{{\cal A}}
\newcommand{\cb}{{\cal B}}
\newcommand{\cc}{{\cal C}}
\newcommand{\cd}{{\cal D}}
\newcommand{\ce}{{\cal E}}
\newcommand{\cg}{{\cal G}}
\newcommand{\ci}{{\cal I}}
\newcommand{\cn}{{\cal N}}
\newcommand{\cs}{{\cal S}}
\newcommand{\ct}{{\cal T}}
\newcommand{\rmd}{\mbox{\rm d}}
\newcommand{\io}{\iota}
\newcommand{\bnot}{\begin{thr} {\bf Notation }}
\newcommand{\lgl}{\langle}
\newcommand{\rgl}{\rangle}
\newcommand{\spp}{\mbox{\rm supp\,}}
\newcommand{\id}{\mathop{\mathrm{id}}}
\newcommand{\pro}{\mathop{\mathrm{pr}}}
\newcommand{\dist}{\mathop{\mathrm{dist}}}
\newcommand{\clb}{\overline{B}_}
\newcommand{\sgn}{\mathop{\mathrm{sgn}}}

\parskip=2mm

\begin{center}
{\bf \Large On the foundations of nonlinear generalized functions II}
\vskip2mm
{\large         M. Grosser}

        Universit\"at Wien\\ 
        Institut f\"ur Mathematik
\end{center}

{\small
{\sc Abstract.} 
This paper gives a comprehensive analysis of algebras of
Colombeau-type generalized functions in the range between the 
diffeomorphism-invariant quotient algebra $\mathcal{G}^d =
\mathcal{E}_M/\mathcal{N}$ introduced in part I
and Colombeau's original algebra $\mathcal{G}^e$. Three main results
are established: First, a simple criterion describing membership in
$\mathcal{N}$ (applicable to all types of Colombeau algebras) is given.
Second, two counterexamples demonstrate that $\mathcal{G}^d$ is not
injectively included in $\mathcal{G}^e$. Finally, it is shown that in
the range  ``between'' $\mathcal{G}^d$ and $\mathcal{G}^e$ only one more
construction leads to a diffeomorphism invariant algebra. 
In analyzing the latter, several classification results essential
for obtaining an intrinsic description of $\mathcal{G}^d$ on manifolds
are derived.

2000 {\it Mathematics Subject Classification}.\hspace{-.5pt} Primary 46F30; 
Secondary 26E15, 46E50, 35D05.

{\it Key words and phrases}. Algebras of generalized functions, Colombeau algebras,
calculus on infinite dimensional spaces, convenient vector spaces, diffeomorphism
invariance.
}

\pagenumbering{arabic}
\setcounter{page}{1}
\setcounter{section}{11}

\section{Introduction to part II}\lb{introII}
In the present article which is the second in a series of two,
we continue the study of diffeomorphism invariant Colombeau
algebras. We will use freely notation and results from the first
part (\cite{fo}); the latter will be referred to herein simply as
``Part~I''. Also, numbering of sections, theorems and formulas will
be continued.

The main result of section \rf{cond0}
permits one to simplify the definition
of the ideal $\cn$ considerably: It dispenses with taking into
account the derivatives of the representative being tested.
This applies to virtually all versions of Colombeau algebras.
This seemingly technical modification, however, has decisive effects on
applications: For example, it makes it considerably easier to
prove uniqueness of the solutions of many differential equations.
Section \rf{calc2} complements section 
4 (``Calculus'') of Part~I
by certain results needed in section \rf{ex}. In particular, it is
shown that $\cc^\infty(U,F)$ is complete with respect to the
topology of uniform convergence (on a suitable family of bounded
sets) in all derivatives resp.\
differentials, provided $F$ is complete as a locally convex space.
In section \rf{ex} we 
show that
the diffeomorphism invariant algebra $\cg^d(\Om)$
presented in section 7 of Part I
is not injectively included in
the Colombeau algebra $\cg^e(\Om)$ of \cite{c2} (which, to be sure, is the standard
version among those being independent from the choice of a particular
approximation of the delta distribution)
by constructing two counterexamples.
In section \rf{spec} we develop a framework allowing to classify
the range of algebras 
which can be positioned between $\cg^d(\Om)$ and
(the smooth version of) $\cg^e(\Om)$.
In particular, we are going to discuss to which extent at least
the definition of the  algebra
introduced by J.~F.~Colombeau and A.~Meril in
\cite{CM}
has to be modified to obtain diffeomorphism invariance.
This leads to the construction of the (diffeomorphism
invariant) Colombeau algebra  $\cg^2(\Om)$  which is closer to the
algebra  of \cite{CM}
than the algebra $\cg^d(\Om)$ (section \rf{g2}). Certain classification
results of sections
\rf{spec} and \rf{g2} are essential for obtaining an intrinsic
description of Colombeau algebras on manifolds (see \cite{vi}).
The concluding section \rf{conclu}
points out that also weaker invariance properties than with respect
to all diffeomorphisms should be envisaged for Colombeau algebras,
in particular regarding applications.

In the following, we will abbreviate $R\circ S^{(\eps)}$ as
$R_\eps$, throughout. Terms of the form $\pa^\al\rmd_1^kR_\eps$
always are to be read as $\pa^\al\rmd_1^k(R_\eps)$.

\section{A simple condition equivalent to negligibility}\lb{cond0}
The property of a representative $R\in\ce(\Om)$
of a generalized function
$[R]\in\cg(\Om)$ to belong to the ideal $\cn(\Om)$ was defined in
7.3. Theorem 18\,(2$^\circ$) of \cite{JEL}
resp.\  Theorem 7.13
of section 7 give an equivalent condition
replacing the term $\pa^\al(R(S_\eps\phi(\eps,x),x))$ occurring in
7.3 by 
$(\pa^\al\rmd_1^k R_\eps)(\vphi,x)(\psi_1,\dots,\psi_k)$.
Moreover, Theorem 18\,(1$^\circ$) of \cite{JEL} shows that
we still get a condition 
equivalent to $R\in\cn(\Om)$ if we
simply omit the differential with respect to the first variable
$\vphi$ from the statement of (2$^\circ$), provided $R$ is assumed to be
moderate. In the following, we are going to show that
a further simplification is possible which might seem rather
drastic at first glance: It is not even necessary to
consider partial derivatives with respect to $x\in\Om$.
In order to facilitate comparing the conditions mentioned 
so far we include all of them in the following theorem,
though only $(0^\circ)$ is new.
\bt\lb{thcond0}
Let $\Om$ be an open subset of $\R^s$ and $R\in\ce_M(\Om)$.
Then each of the following conditions is equivalent
to $R\in\cn(\Om)$ (in the sense of 7.3):\ms
$(0^\circ)$\quad
$\forall K\subset\subset\Om\  
   \forall n\in\N\ \exists q\in\N
   \ \forall B\,(\mbox{bounded})\,\subseteq \cd(\R^s)$:
$$R_\eps(\vphi,x)=O(\eps^n)
   \qquad\qquad (\eps\to0)$$
uniformly for $x\in K$, $\vphi\in B\cap\ca_q(\R^s)$.\ms
$(1^\circ)$\quad
$\forall K\subset\subset\Om\ \forall\al\in\N_0^d\ 
   \forall n\in\N\ \exists q\in\N
   \ \forall B\,(\mbox{bounded})\,\subseteq \cd(\R^s)$:
$$\pa^\al R_\eps(\vphi,x)=O(\eps^n)
   \qquad\qquad (\eps\to0)$$
uniformly for $x\in K$, $\vphi\in B\cap\ca_q(\R^s)$.\ms
$(2^\circ)$\quad
$\forall K\subset\subset\Om\ \forall\al\in\N_0^d\ 
   \forall k\in\N_0\ \forall n\in\N\ \exists q\in\N
   \ \forall B\,(\mbox{bounded})\,\subseteq \cd(\R^s)$:
$$\pa^\al\rmd_1^k R_\eps(\vphi,x)(\psi_1,\dots,\psi_k)=O(\eps^n)
   \qquad\qquad (\eps\to0)$$
uniformly for $x\in K$, $\vphi\in B\cap\ca_q(\R^s)$,
   $\psi_1,\dots,\psi_k\in B\cap\ca_{q0}(\R^s)$.
\et
\pr
The equivalence of each of $(1^\circ)$ and $(2^\circ)$
with $R\in\cn(\Om)$ is a part of
Theorem 18 of \cite{JEL}. $(1^\circ)\Rightarrow(0^\circ)$ 
being trivial, it remains to show
$(0^\circ)\Rightarrow(1^\circ)$.
To this end, we will prove, assuming $R\in\ce_M(\Om)$ to satisfy
$(0^\circ)$, that $R$ satisfies $(1^\circ)$ for
$\al:=e_i$, i.e., $\pa^\al=\pa_i$ ($i=1,\dots,s$)
and that, in addition, $\pa_i R$ again is moderate and satisfies 
$(0^\circ)$. Then it will follow by induction that $(1^\circ)$
holds for all $\al\in\N_0^s$.

So suppose $R\in\ce_M(\Om)$ to satisfy $(0^\circ)$
and let $K\subset\subset\Om$ and
$n\in\N$ be given. For $\de:=\min(1,\dist(K,\pa\Om))$,
set $L:=K+\clb{\frac{\de}{2}}(0)$.
Then $K\subset\subset L\subset\subset\Om$.
Now by moderateness of $R$ and Theorem 7.12, choose $N\in\N$
such that for every bounded subset $B$ of $\cd(\R^s)$ the relation
$\pa_i^2R_\eps(\vphi,x)=O(\eps^{-N})$ as $\eps\to0$
holds,
uniformly for $x\in L$, $\vphi\in B\cap\ca_0(\R^s)$. Next,
by the assumption of $(0^\circ)$ to hold for $R$, choose
$q\in\N$ such that, again for every bounded subset $B$ of
$\cd(\R^s)$, we have
$R_\eps(\vphi,x)=O(\eps^{2n+N})$ as $\eps\to0$,
uniformly for $x\in L$, $\vphi\in B\cap\ca_q(\R^s)$.
Now suppose a bounded subset $B$ of $\cd(\R^s)$ to be given;
let $\vphi\in B\cap\ca_q(\R^s)$, $x\in K$ and
$0<\eps<\frac{\de}{2}$; hence $x+\eps^{n+N}e_i\in L$.
By Taylor's theorem, we conclude
(to be precise, separately for the real and imaginary part of $R$)
\beas
R_\eps(\vphi,x+\eps^{n+N}e_i)
   &=&R_\eps(\vphi,x)+\pa_iR_\eps(\vphi,x)\eps^{n+N}+
      \frac{1}{2}\pa_i^2R_\eps(\vphi,x_\theta)\eps^{2n+2N}
\eeas
where $x_\theta=x+\theta\eps^{n+N}e_i$ for some $\theta\in(0,1)$;
note that also $x_\theta\in L$.
Consequently,
\beas
\pa_iR_\eps(\vphi,x)
   &=&
       \underbrace
         {\left(R_\eps(\vphi,x+\eps^{n+N}e_i)
          -R_\eps(\vphi,x)\right)}
         _{O(\eps^{2n+N})} \eps^{-n-N}
       -\underbrace
         {\frac{1}{2}\pa_i^2R_\eps(\vphi,x_\theta)}
         _{O(\eps^{-N})}
       \eps^{n+N},
       \eeas
uniformly for $\vphi\in B\cap\ca_q(\R^s)$, $x\in K$.
Having demonstrated $\pa_iR_\eps(\vphi,x)=O(\eps^n)$
for all $i=1,\dots,s$, observe that
$\pa_i(R_\eps)=(\pa_iR)_\eps$. Therefore,
$\pa_iR$ again satisfies $(0^\circ)$.
According to Theorem 7.10 (which is non-trivial,
see the discussion in section 7),
$\pa_iR$ is also moderate .
By the remark made above, this
completes the proof.
\ep

The reader acquainted with E.~Landau's paper \cite{landau} will
easily recognize the method employed therein to form the basis of
the preceding proof (though not mentioned explicitly in \cite{JEL},
this equally applies to the proof of $(1^\circ)\Rightarrow(2^\circ)$
of Theorem 18 of \cite{JEL}).

The part of Theorem \ref{thcond0}
saying that for moderate functions (the appropriate analog of)
condition $(0^\circ)$ is equivalent to 
negligibility
applies to virtually all versions of
Colombeau algebras of practical importance,
in particular, to the following:
\begin{itemize}
\item
For the special algebra as defined,
e.g., in \cite{MObook}, p.~109, just replace the term
$R_\eps(\vphi,x)$
in condition $(0^\circ)$ by $u_\eps(x)$.
\item
For the classical full Colombeau algebra
of \cite{c2} simply drop the uniformity requirement
concerning $\vphi$ from $(0^\circ)$.
\item
For the diffeomorphism invariant Colombeau algebra
$\cg^2(\Om)$ to be introduced in section \rf{g2},
the corresponding result is stated as Theorem \rf{thcond0g2}.
\item
For the special algebra on smooth mainfolds
the corresponding result follows from the local
characterization of generalized functions
(see \cite{roldiss}, 4.2).
\item
The latter also applies to the intrinsically defined
full Colombeau algebra on manifolds (\cite{vi}, Corollary 4.5).
\end{itemize}
In the first and second of these four instances, the
respective proofs are obtained by appropriately slimming down
the corresponding argument of the proof of Theorem
\rf{thcond0}.

The seemingly technical difference
between $(0^\circ)$ and the remaining conditions
(including negligibility of $R$) has decisive effects on
applications: For example, if the uniqueness of a
solution of a differential equation is to be shown one supposes
$R_1,R_2$ to be representatives of solutions. Note that this
includes the assumption that $R_1,R_2\in\ce_M(\Om)$, hence
Theorem \rf{thcond0} may be applied. For $[R_1]=[R_2]$ in $\cg(\Om)$
we have to show that $R:=R_1-R_2\in\cn(\Om)$. Now it suffices to check
condition $(0^\circ)$ rather than $(1^\circ)$ (resp.\ $(2^\circ)$
resp.\ the original definition of $R\in\cn(\Om)$),
i.e., there is no need
to analyze the behaviour of any derivative of $R$.

Apart from that, condition $(0^\circ)$ is also of theoretical relevance.
To give a sample, let us demonstrate
that it allows to simplify considerably the proof of
statement (iv) of Theorem 7.4 in Part I
(saying that $(\io-\si)(\cc^\infty(\Om))\subseteq\cn(\Om)$):
Since
$\io f-\si f\in\ce_M(\Om)$ by (i) and (ii) of Theorem 7.4,
it is sufficient for $\io f-\si f\in\cn(\Om)$
to show that
$$(\io f-\si f)(S_\eps\vphi,x)=\int\limits_{\frac{\Om-x}{\eps}}
  \left[f(z\eps+x)-f(x)\right]\vphi(z)\,dz
  =O(\eps^{q+1}),$$
uniformly for $x\in K$ and $\vphi$ ranging over some bounded subset
of $\ca_q(\R^d)$. This, however, is immediate.

\section{Some more calculus}\lb{calc2}
Both the counterexamples to be constructed in section \rf{ex}
will take the
form of infinite series, being absolutely convergent in each
derivative. Thus we need a theorem
guaranteeing the completeness of
$\ce(\Om)=\cc^\infty(U(\Om),\C)\equiv
    \cc^\infty(\ca_0(\Om)\times\Om,\C)$
with respect to the corresponding topology.
The remarkable
ease of the proof of this generalization of a standard result of
elementary real analysis clearly exhibits the virtues of calculus
in convenient vector spaces as outlined in section 4.
To this end, let $E,F$ be locally convex spaces
and $U$ an open subset of $E$. If $f:U\to F$ is
smooth, its
\linebreak
$n$-th differential $\rmd^n\!f$ belongs to
$\cc^\infty(U,L^n(E^n,F))$ where $L^n(E^n,F)$ denotes the space
$L(E,\dots,E;F)$ of $n$-linear bounded maps from
$E\times\dots\times E$ ($n$ factors) into $F$.
(For $n=0$, set $L^n(E^n,F):=F$.)
On $\cc^\infty(U,L^n(E^n,F))$, let $\tau_{cb}^n$ denote
the topology of uniform $F$-convergence on subsets of the form $K\times B$
where $K$ is a compact subset of $U$ and $B$ is bounded in
$E^n=E\times\dots\times E$. Let $\cc^\infty(U,F)$ carry the initial
(locally convex) topology $\tau^\infty$ induced by the family
$(\rmd^n,\cc^\infty(U,L^n(E^n,F)),\tau_{cb}^n)_{n\ge0}$, i.e., the
topology of uniform convergence of all derivatives
(that is to say, differentials) on sets $K\times B$ as above.
Note that on $\cc^\infty(\R,F)$, $\tau^\infty$ 
is just the usual Fr\'echet topology of compact convergence in all
derivatives.

\bt\lb{cinftcomp}
Let $E,F$ be locally convex spaces, assume $F$ to be complete
and let $U$ be an open subset of $E$. Then $\cc^\infty(U,F)$
is complete with respect to the topology
$\tau^\infty$ of uniform $F$-convergence of all differentials
on subsets of the form $K\times B$
where $K$ is a compact subset of $U$ and $B$ is bounded in
the appropriate product $E^n=E\times\dots\times E$. Moreover,
for each $p\in\N$, the operator
$\rmd^p:\cc^\infty(U,F)\to\cc^\infty(U,L^p(E^p,F))$ is continuous
if both the domain and the range space carry
the respective topology $\tau^\infty$.
\et
\pr
Let $(f_\io)$ be a net in $\cc^\infty(U,F)$ which is Cauchy
with respect to $\tau^\infty$, that is, suppose $(\rmd^n\!f_\io)$
to be a Cauchy net in $\cc^\infty(U,L^n(E^n,F))$ with respect
to $\tau^n_{cb}$ for each $n=0,1,2,\dots$.
Due to the completeness of $F$, each net $(\rmd^n\!f_\io)$
has a limit $f^{[n]}:U\times E^n\to F$
with respect to (the obvious extension of) $\tau^n_{cb}$.
In particular, $(f_\io)$ converges to
some function $f:=f^{[0]}:U\to F$.
Consider a smooth curve $c:\R\to U$;
then for each $\io$, $f_\io\circ c$ is smooth from $\R$ to $F$,
its $n$-th derivative at $t\in\R$
being given as a certain sum of terms of the
form $(\rmd^l\!f_\io)(c(t))(c^{(k_1)}(t),\dots,c^{(k_l)}(t))$
where $1\le l\le n$ and $\sum k_j=n$, due to the
chain rule. With $t$ ranging over some compact subset of $\R$,
the values attained by $c^{(k)}(t)$ form a compact
subset of $U$ resp.\ $E$, for each $k\in\N_0$.
Now it follows from the Cauchy property of $(f_\io)$
that $(f_\io\circ c)$ is Cauchy in $\cc^\infty(\R,F)$
with respect to uniform convergence
in all derivatives on compact sets.
From the completeness of the latter space we conclude that
the limit of $(f_\io\circ c)$
exists as a smooth function and is equal to $f\circ c$.
This argument being valid for any smooth curve $c$,
$f$ itself is smooth. To establish $f=\lim f_\io$ with respect to
$\tau^\infty$,
it remains to show that for any $n\in\N$,
$\rmd^n\!f=f^{[n]}$, i.e., that for all $x\in U$,
$v_1,\dots,v_n\in E$ we have
\be\lb{dnf}
(\rmd^n\!f)(x)(v_1,\dots,v_n)=\lim\,(\rmd^n\!f_\io)(x)(v_1,\dots,v_n),
\end{equation}
For a straight line $c(t)=x+tv$ we obtain, at $t=0$,
$(g\circ c)^{(n)}(0)=(\rmd^n\!g)(x)(v,\dots,v)$
for any $g\in\cc^\infty(U,F)$. Therefore,
$$(\rmd^n\!f)(x)(v,\dots,v)=(f\circ c)^{(n)}(0)
   =\lim\,(f_\io\circ c)^{(n)}(0)
   =\lim\,(\rmd^n\!f_\io)(x)(v,\dots,v).$$
Equation (\rf{dnf}) now follows by polarization
(see, e.g., \cite{KM}, Lemma 7.13\,(1)).
Finally, the continuity of $\rmd^p$ with respect to the initial
topologies $\tau^\infty$ is immediate from the following
commutative diagram:
\[
\begin{CD}
(\cc^\infty(U,F),\tau^\infty)
   @>{\rm d}^p>>(\cc^\infty(U,L^p(E^p,F)),\tau^\infty)\\
@VV{\rm d}^{p+n}V       @VV{\rm d}^nV                   \\
(\cc^\infty(U,L^{p+n}(E^{p+n},F)),\tau^{p+n}_{cb})
   @>\id>>(\cc^\infty(U,L^n(E^n,L^p(E^p,F))),\tau^n_{cb}) \\
\end{CD}
\]
Observe that the lower horizontal arrow is a linear
homeomorphism, due to $L^p(E^p,F)$ carrying the topology of uniform
convergence on bounded sets.
\ep

For the rest of this section, let $U$ denote a (non-empty) open subset
of a closed affine subspace $E_1$ of some locally convex space $E$,
$E_0$ the linear subspace parallel to $E_1$ and $F$ a complete
locally convex space. {\it Mutatis mutandis}, Theorem
\rf{cinftcomp} is valid also in this slightly more general
situation. The vectors
$v_1,\dots,v_n$ to be plugged into $\rmd^n\!f(x)$ now have to be
taken from $E_0$, as well as $B$ has to denote a bounded subset
of $E_0^n$.

\bd\lb{defeb} Assume, in addition to the above, that the topology
of $F$  is generated by some family ${\cal P}$ of semi-norms.
For fixed $n\in\N_0$, let $(f_k)_{k\in\N}$ denote a sequence of
functions
\beas
f_k:U\times E_0\times \dots \times E_0 &\to& F  \hspace{100pt}
              (n\ \mbox{\rm factors}\ E_0)\\
f_k:(\ x\ ,\ v_1\ ,\dots,\ v_n\ )&\mapsto&f_k(x)(v_1,\dots,v_n).
\eeas
We say that $(f_k)$ is \mbox{\rm exponentially bounded on}
$U\times E_0^n$  (an
\mbox{\rm(eb)}-sequence, for short) if for each compact subset $K$ of $U$,
each bounded subset $B$ of $E_0$ and each $p\in{\cal P}$
there exists a constant $C(\ge1)$ such that
$p(f_k(x)(v_1,\dots,v_n))\le C^k$
for any $k\in\N$, $x\in K$ and $v_i\in B$ $(i=1,\dots,n)$.
\et
Define $(f_k)+(g_k):=(f_k+g_k)$ and $\lambda(f_k):=(\lambda f_k)$
($\lambda\in\C$), as well as \linebreak $(f_k)\cdot(g_k):=(f_k\cdot g_k)$
provided $F$ is a (complete) locally convex topological algebra.
Then the following proposition is immediate, due to
$C_1^k+C_2^k\le(C_1+C_2)^k$ and $C_1^kC_2^k=(C_1C_2)^k$:
\bp\lb{ebalg}
The set of \mbox{\rm(eb)}-sequences forms a linear space
(resp.\ an algebra if $F$ is a locally convex
algebra with jointly sequentially continuous multiplication)
with respect to the
operations defined above.
\et
\bt\lb{theb}
Let $f_k\in\cc^\infty(U,F)$ for every $k\in\N$. Assume that for
each fixed $n\in\N_0$, $(\rmd^n\!f_k)_k$ is \mbox{\rm(eb)}
on $U\times E_0^n$. Then
$\sum\limits_{k=0}^\infty\frac{1}{k!}f_k$ is convergent with
respect to $\tau^\infty$ to some $f\in\cc^\infty(U,F)$. Moreover,
$\rmd^n\!f=\sum\limits_{k=0}^\infty\frac{1}{k!}\,\rmd^n\!f_k$ for
every $n\in\N_0$ where also the latter series converges with
respect to $\tau^\infty$.
\et
\pr
Fix $n\in\N_0$,
a compact subset $K$ of $U$ and a bounded subset $B$ of $E_0$.
Since $(\rmd^n\!f_k)_k$ is (eb),
$\sum\limits_{k}\frac{1}{k!}\,\rmd^n\!f_k$ is majorized,
uniformly on $K\times
B^n$, by $\sum\limits_{k}\frac{C_n^k}{k!}$ for some constant
$C_n(\ge1)$ depending only on $n$, $K$ and $B$. Consequently,
$\sum\limits_{k}\frac{1}{k!}f_k$ is Cauchy with respect to
$\tau^\infty$. Now both the convergence of
$\sum\limits_{k}\frac{1}{k!}f_k$ and the admissibility of term-wise
differentiation follow from Theorem \rf{cinftcomp}.
\ep

In the sequel, \rf{defeb}--\rf{theb} will only be used for
$F=\C$; the extension to locally convex algebras being for
free virtually, we chose to state them in the general form
to indicate the scope of Theorem \rf{theb}.

\section{Non-injectivity of the canonical homomorphism from
$\cg^d(\Om)$ into $\cg^e(\Om)$}\lb{ex}

For every open subset $\Om$ of $\R^s$, there is a canonical
algebra homomorphism $\Phi$ from the diffeomorphism invariant
Colombeau algebra $\cg^d(\Om)$ of \cite{JEL} (see section
7) to the ``classical'' (full) Colombeau algebra
$\cg^e(\Om)$ introduced in \cite{c2}, 1.2.2 (the upper subscript $e$
being taken from the title ``Elementary Introduction to New
Generalized Functions'' of the latter monograph). In this section,
we are going to show that $\Phi$ is not injective in general by
constructing a representative $R$ of a generalized function
$[R]\in\cg^d(\Om)$ such that $[R]\neq0$, yet $\Phi[R]=0$.

By superscripts $d,e$ we will distinguish between ingredients (as
listed in section 3) for constructing 
$\cg^d(\Om)$ resp.\ $\cg^e(\Om)$. Observe that superscripts
$d,e$ are independent of superscripts $J,C$ as introduced in
section 5: Each of the (non-isomorphic) algebras
$\cg^d(\Om)$, $\cg^e(\Om)$ has equivalent descriptions
in the C- and the J-formalism, respectively.
As in section 7 of Part I, we will use the C-formalism
also in the present context.
All the relevant definitions are
to be found in section 7 (for $\cg^d(\Om)$) resp.\
\cite{c2} (for $\cg^e(\Om)$). For the present purpose, the
following of them are of particular importance:
\bea U^d(\Om)&:=&T^{-1}(\ca_0(\Om)\times\Om)\nonumber\\
     U^e(\Om)&:=&T^{-1}(\ca_1(\Om)\times\Om)\footnotemark\lb{ueom}
\eea
\beas \ce^d(\Om)&:=&\cc^\infty(U^d(\Om))\\
     \ce^e(\Om)&:=&\{R:U^e(\Om)\to\C\mid x\mapsto R(\vphi,x)
     \ \mbox{\rm is smooth on $U_\vphi$ for each }\vphi\}
\eeas
where $U_\vphi$
denotes the (open) set $\{x\mid (\vphi,x)\in U^e(\Om)\}$.
\footnotetext{The choice of $\ca_1(\Om)$ rather than $\ca_0(\Om)$
in the definition of $U^e(\Om)$ is due to Co\-lom\-beau
(\cite{c2}, 1.2.1). We decided to keep the original form of
$\cg^e(\Om)$ although all the results of this section would remain
valid (and, in fact, even slightly easier to formulate)
choosing also $U^e(\Om)$ to be $T^{-1}(\ca_0(\Om)\times\Om)$.}

From now on, we will omit
specifying $\Om$ explicitly whenever it is clear which domain is
intended. Let $j:U^e\to U^d$ denote set-theoretic inclusion.
To see that the restriction map $\Phi_0=j^*$ maps $\ce^d$ into
$\ce^e$ we have to pass from C-representatives to
J-representatives: Smoothness of $R^d\in\ce^d$, by definition,
is equivalent to smoothness of
$(T^*)^{-1}R^d\in\cc^\infty(\ca_0(\Om)\times\Om)$ while
for $R^e\in\ce^e$, smoothness of $x\mapsto R^e(\vphi,x)$
is equivalent to smoothness of
$x\mapsto(T^*)^{-1}R^e(\vphi(.-x),x)$. From this it is clear that
$\Phi_0R^d\in\ce^e$ for $R^d\in\ce^d$.

$\Phi_0$ even maps $\ce^d_M$ into $\ce^e_M$ and
$\cn^d$ into $\cn^e$, respectively. This follows easily by
inspecting the corresponding definitions: For $R^d\in\ce^d$,
$R^e\in\ce^e$; we have, by definition (omitting the quantifiers
``$\forall K\subset\subset\Om\ \forall\al\in\N_0^s\ \exists
N\in\N$''),
$$\ba{rlll}
R^d\in\ce_M^d&\Leftrightarrow&
   \forall\phi\in\cc_b^\infty(I \times\Om,\ca_0(\R^s)):&
   \sup\limits_{x\in K}|\pa^\al(R^d(S_\eps\phi(\eps,x),x))|
   =O(\eps^{-N})\\
R^e\in\ce_M^e&\Leftrightarrow&
   \forall\vphi\in\ca_N(\R^s):&
   \sup\limits_{x\in K}|\pa^\al(R^e(S_\eps\vphi,x))|
   \hphantom{(\eps,x)}=O(\eps^{-N})
\ea$$
Obviously, each test object $\vphi\in\ca_N(\R^s)$ can be viewed as
a particular case of a test object
$\phi\in\cc_b^\infty(I \times\Om,\ca_0(\R^s))$
by setting $\phi(\eps,x):=\vphi$ independently of $\eps,x$.
Thus from $R^d\in\ce_M^d$ it follows that $\Phi_0R^d\in\ce_M^e$.
A similar argument shows that $\Phi_0R^d\in\cn^e$ provided
$R^d\in\cn^d$. Note that the condition for the membership
of $R^e$ in $\cn^e$
as given in \cite{c2}, 1.1.11, that is
(this time omitting
``$\forall K\subset\subset\Om\ \forall\al\in\N_0^s$'')
$$
   \exists N\ \exists\ga:\N\to\R\ \forall q\ge N
   \ \forall\vphi\in\ca_q(\R^s):
   \sup\limits_{x\in K}|\pa^\al(R^e(S_\eps\vphi,x))|
   =O(\eps^{\ga(q)-N})
$$
(where $\ga(q)\nearrow\infty$) is easily seen to be equivalent to
$$
\forall n\ \exists q\ 
   \hphantom{mmmmmmmmm}
   \forall\vphi\in\ca_q(\R^s):
   \sup\limits_{x\in K}|\pa^\al(R^e(S_\eps\vphi,x))|
   =O(\eps^n)\hphantom{\ga(q)}
$$
which has the same structure as the condition in 7.3
for $R^d$ to belong to $\cn^d$:
$$
\hspace*{0pt}\hphantom{m}\forall n\ \exists q
  \ \hphantom{mii}\forall\phi\in\cc_b^\infty(I \times\Om,\ca_q(\R^s)):
   \sup\limits_{x\in K}|\pa^\al(R^d(S_\eps\phi(\eps,x),x))|
   =O(\eps^{n})\hphantom{\ }.
$$
Due to the invariance of $\ce_M$ and $\cn$ under $\Phi_0$, $\Phi_0$
induces a map $\Phi:\cg^d(\Om)\to\cg^e(\Om)$ acting on representatives
as restriction from $T^{-1}(\ca_0(\Om)\times\Om)$ to
$T^{-1}(\ca_1(\Om)\times\Om)$. $\Phi$ is an algebra homomorphism
respecting the embeddings of $\cd'(\Om)$ and differentiation.

\brem (i) If we had chosen to set
$U^e(\Om)=T^{-1}(\ca_0(\Om)\times\Om)$ (contrary to \cite{c2},
cf. the footnote to (\rf{ueom}) above) $j$ would be the identity
map on $U^d(\Om)=U^e(\Om)$ and $\Phi_0$ would be set-theoretic
inclusion, hence injective.\smallskip\\
(ii) Regarding the question of injectivity of $\Phi$, the fact that
$\ca_1(\Om)$ has been used in \cite{c2} and in (\rf{ueom}) above
to define $U^e(\Om)$ (as compared to $\ca_0(\Om)$ in \cite{JEL}
for defining $U^d(\Om)$) is completely irrelevant: Although this
choice renders $\Phi_0$ non-injective  in general (consider
$(0\neq)R\in\cc^\infty(U^d(\R))=\cc^\infty(\ca_0(\R)\times\R)$
given by $(\vphi,x)\mapsto\int\xi\vphi(\xi)\,d\xi$:
$\Phi_0R=0$ by the very definition of $U^e(\R)=\ca_1(\R)\times\R$),
${\cal M}:=\ker\Phi_0$ is contained in $\cn^d$ since each
$R\in\ker\Phi_0$ vanishes identically on pairs $(\phi(\eps,x),x)$
where $\phi$ is a test object taking values in $\ca_q(\R)$
($q\ge1$). Thus the canonical image of ${\cal M}$ in
$\cg^d:=\ce_M^d\big/\cn^d$ is trivial.
\et

Having discussed $\Phi$ in detail, we will omit $j$ and $\Phi_0$
from our notation in the sequel.
Now we can state precisely which properties a function
$R:U^d(\Om)\to\C$ has to satisfy if it is to refute the injectivity
of $\Phi$:
\begin{itemize}
\item[(i)] $R\in\ce^d$, i.e., $R$ has to be smooth;
\item[(ii)] $R\in\ce_M^d$,
\item[(iii)]$R\notin\cn^d$,
\item[(iv)] $R\in\cn^e$.
\end{itemize}
In the following, we will define maps $P,Q:U^d(\R)\to\C$ each of
which satisfies (i)--(iv) above, thereby providing a counterexample
to the conjecture of the canonical map $\Phi$ being injective. We
will give the complete argument for $P$ while only indicating how
to adapt the proof to get the analogous result for $Q$.

For the definition of $P,Q$ let $s:=1$, $\Om:=\R$. We continue using
the C-formalism. Although now
$U(\Om)=\ca_0(\R)\times\R=\ca_0(\Om)\times\Om$ note that the
C-formalism, nevertheless, differs from the J-formalism with
respect to embedding $\cd'$, differentiation, testing (which
involves $T$ in the case of the J-formalism) and, finally, with
respect to the action induced by a diffeomorphism. As a prerequisite
for writing down $P,Q$ explicitly, we introduce the following
notation:
$$\ba{lrcll}
&\lgl\vphi|\vphi\rgl&:=&\int\vphi(\xi)\overline{\vphi(\xi)}\,d\xi
      \qquad\qquad &(\vphi\in\cd(\R))\medskip\\
v_k\in\cd'(\R):\qquad&\lgl v_k,\vphi\rgl&:=&\int\xi^k\vphi(\xi)\,d\xi
      \qquad\qquad&(\vphi\in\cd(\R),\ k\in\N_0)\medskip\\
v_\frac{1}{2}\kern-2pt\in\cd'(\R):\qquad&\lgl v_\frac{1}{2},
                     \vphi\rgl&:=&\int|\xi|^\frac{1}{2}\vphi(\xi)\,d\xi
      \qquad\qquad&(\vphi\in\cd(\R))\medskip\\
&v(\vphi)&:=&\lgl\vphi|\vphi\rgl^\frac{1}{2}
                              \lgl v_\frac{1}{2},\vphi\rgl
      \qquad\qquad&(\vphi\in\cd(\R))\medskip\\
&g(x)&:=&\frac{x}{1+x^2}&(x\in\R)\medskip\\
&e(x)&:=&\cases 
             \exp(-\frac{1}{x})\quad&(x>0)\\
             0                     &(x\le0)
          \endcases &(x\in\R)\medskip\\
&\ga_k&:=&k+\frac{1}{k}&(k\in\N).
\ea$$
Finally, choose an (even) function
$\si\in\cd(\R)$ satisfying $0\le\si\le1$,
$\si(x)\equiv1$ for $|x|\le\frac{1}{2}$,
$\si(x)\equiv0$ for $|x|\ge\frac{3}{2}$ and set
$$h_k(x):=\si(x)\cdot 2g(x)+(1-\si(x))\cdot
\sgn(x)\cdot|2g(x)|^{\ga_k}\qquad\qquad(x\in\R,
 \ k\in\N).$$
Being bounded and linear resp.\ bilinear (over $\R$),
$v_k$, $v_\frac{1}{2}$ and $\lgl\,.\,|\,.\,\rgl$ are smooth on $\cd(\R)$ ($k\in\N_0$).
On $\ca_0(\R)$, $\lgl\vphi|\vphi\rgl>0$. Thus $v$ is smooth on
$\ca_0(\R)$ as a product of smooth functions. Observe that
$\ca_q(\R)=\ca_0(\R)\cap\bigcap\limits_{k=1}^{q}\ker v_k$.

In the sequel, we will make use of the following facts
concerning $g$ and $e$:
For every $n\in\N_0$ there exists a constant $c_n>0$ such that
for all $x\ne0$
$$|g^{(n)}(x)|\le\frac{c_n}{|x|^{n+1}}.$$ 
The derivatives of $e$ can be written in the following form:
$$e^{(n)}(x)=e(x)\cdot\frac{q_n(x)}{x^{2r}}=\cases 
             \exp(-\frac{1}{x})\cdot\frac{q_n(x)}{x^{2r}}\quad&(x>0)\\
             0                     &(x\le0)
          \endcases$$ 
for every $n\in\N$ where $q_n$ is a polynomial of degree $n-1$ and
$\frac{0}{0}:=0$. 

Scaling of $\vphi$ produces the following relations:
$$\ba{rcl} 
     \lgl S_\eps\vphi|S_\eps\vphi\rgl&=&
                 \frac{1}{\eps}\lgl\vphi|\vphi\rgl\medskip\\
     \lgl v_k,S_\eps\vphi\rgl&=&
                  \eps^k\lgl v_k,\vphi\rgl\medskip\\
     \lgl v_\frac{1}{2},S_\eps\vphi\rgl&=&
                  \eps^{\frac{1}{2}}\lgl v_\frac{1}{2},\vphi\rgl
                                      \medskip\\
     v(S_\eps\vphi)&=&
                  v(\vphi).
\ea$$
Apart from abbreviating $R\circ S^{(\eps)}=R\circ(S_\eps\times\id)$
as $R_\eps$ for any function $R$ defined on $\ca_0(\R)\times\R$,
we also will write $R_\eps$ for $R\circ S_\eps$ if $R$ is defined on
$\ca_0(\R)$.
\bd \lb{defpq}
Let $\vphi\in\ca_0(\R)$, $x\in\R$ and set
\bea
\lb{defp}
P(\vphi,x)&:=&\sum_{k=1}^{\infty}\frac{1}{k!}\cdot
   g\big(\lgl\vphi|\vphi\rgl^{\ga_k} e(v(\vphi))\big)\cdot
   \lgl\vphi|\vphi\rgl^{\ga_k}\cdot
   \lgl v_k,\vphi\rgl,\\
\lb{defq}
Q(\vphi,x)&:=&\sum_{k=1}^{\infty}\frac{1}{k!}\cdot
   h_k\big(\lgl\vphi|\vphi\rgl^{\frac{3}{2}}\, \lgl v_{\frac{1}{2}},
                                   \vphi\rgl\big)\cdot
   \lgl\vphi|\vphi\rgl^{\ga_k}\cdot
   \lgl v_k,\vphi\rgl.
\eea
\et
Hence $P$ and $Q$, in fact, only depend on $\vphi$. We will see
below that both series for $P$ and $Q$ converge uniformly on
bounded subsets of $\ca_0(\R)$, making $P$ and $Q$ well-defined.
For $k\in\N$, $\vphi\in\ca_0(\R)$ set
$$P_k(\vphi):=g\big(\lgl\vphi|\vphi\rgl^{\ga_k}
           e(v(\vphi))\big)\cdot
   \lgl\vphi|\vphi\rgl^{\ga_k}\cdot
   \lgl v_k,\vphi\rgl.$$
Fix a positive number $\eta\le1$.
To establish properties (i) and (ii) (i.e., smoothness and
moderateness) of $P$, we will derive estimates of the form
\bea\lb{estpk}
|(\rmd^n(P_k)_\eps)(\vphi)(\psi_1,\dots,\psi_n)|\le
   C_n^k\cdot\eps^{-\frac{1}{k}-n\eta}
\eea
for some constants $C_n\ge1$ not depending on $\eps$
($n\in\N_0$), uniformly on any bounded
subset $B$ of $\cd(\R)$ and $\vphi\in B\cap\ca_0(\R)$,
$\psi_1,\dots,\psi_n\in B\cap\ca_{00}(\R)$. Setting $\eps=1$
in (\rf{estpk})
shows that for each $n\in\N_0$, $(\rmd^n\!P_k)$ is an (eb)-sequence on
$\ca_0(\R)\times\ca_{00}(\R)^n$ which, by Theorem \rf{theb},
implies smoothness of $P$. Considering arbitrary values of
$\eps\in I$, on the other hand, will lead to the proof of
moderateness of $P$.

\subsection{Proof of the estimates (\rf{estpk})}
Fix $n\in\N_0$, $\eps\in I$ and $0<\eta\le1$.
Set
\beas P_k^{(1)}(\vphi)&:=&g\big(\lgl\vphi|\vphi\rgl^{\ga_k}
          e(v(\vphi))\big)\\
      P_k^{(2)}(\vphi)&:=&\lgl\vphi|\vphi\rgl^{\ga_k}\\
      P_k^{(3)}(\vphi)&:=&\lgl v_k,\vphi\rgl.
\eeas
(\rf{estpk}) is equivalent to saying that the functions
$\eps^{\frac{1}{k}+n\eta}\cdot\rmd^n(P_k)_\eps$ form an
(eb)-sequence, with the respective constants in the estimate
independent of $\eps$. In order to prove this,
by Leibniz' rule for the differential of a product
and by Proposition \rf{ebalg} it suffices to 
show that each of the sequences (indexed by $k\in\N$)
$\eps^{n\eta}\rmd^m (P_k^{(1)})_\eps$,
$\eps^{\ga_k}\rmd^m (P_k^{(2)})_\eps=\rmd^m (P_k^{(2)})$ and
$\eps^{-k}   \rmd^m (P_k^{(3)})_\eps=\rmd^m (P_k^{(3)})$
is (eb) for $m\le n$, independently of $\eps$.
In the following, we are going to verify these claims step by
step, starting with the elementary building blocks of the series
defining $P$ resp.\ $Q$.
\smallskip

{\bf Remark.} For $w\in\cd'(\Om)$, $(\rmd w)(\vphi)(\psi)=\lgl w,\psi\rgl$
and $\rmd^l w=0$ for $l\ge2$, due to
the linearity of $w$. $\lgl\,.\,|\,.\,\rgl$ being bilinear over $\R$,
we obtain
$(\rmd\lgl\,.\,|\,.\,\rgl)(\vphi)(\psi)=\lgl\psi|
                                  \vphi\rgl+\lgl\vphi|\psi\rgl$,
$(\rmd^2\lgl\,.\,|\,.\,\rgl)(\vphi)(\psi_1,\psi_2)=\lgl\psi_1|
                                  \psi_2\rgl+\lgl\psi_2|\psi_1\rgl$
and $\rmd^l\lgl\,.\,|\,.\,\rgl=0$ for $l\ge3$.
\bp\lb{listeb}
The following sequences of functions of $\vphi$
(indexed by $k\in\N$) are \mbox{\rm (eb)} (1.\ and 3.\
on $\cd(\R)$, 2.\ on $\ca_0(\R)$):
\ben
\item $\lgl\xi^k,\vphi(\xi)\rgl$,\quad$\lgl|\xi|^k,\vphi(\xi)\rgl$,
      \quad $\lgl\vphi|\vphi\rgl^k$,\quad
      $\lgl\vphi|\vphi\rgl^{\ga_k}$;
\item $\lgl\vphi|\vphi\rgl^{-k}$,\quad 
      $\lgl\vphi|\vphi\rgl^{-\ga_k}$;
\item $(\bet_k)_n:=\bet_k(\bet_k-1)\dots(\bet_k-n+1)$\qquad
               (for fixed $n\in\N_0$; $(\bet_k)_0:=1$)
\een
where the numbers $\bet_k\in\R$ occurring in 3.\ only have to satisfy an
estimate of the form $|\bet_k|\le pk$ for some fixed $p\in\N$. 
\et
\pr
Fix a bounded subset $B$ of $\cd(\R)$ containing at least one
$\vphi\ne0$ . Then there exists a bounded
set $L\subseteq\R$ containing the supports of all $\vphi\in B$.
Let $\mbox{\rm m}(L)>0$ denote the Lebesgue measure of $L$ and
set $C_1:=\max(1,\sup\limits_{\xi\in L}|\xi|)$,
$C_2:=\max (1,\mbox{\rm m}(L))$. Moreover,
$C_3:=\max(1,\sup\limits_{\vphi\in B}\|\varphi\|_\infty)$ is finite.
Now let $\vphi\in B$.

1.\
We have
$$\max\big(|\lgl\xi^k,\vphi(\xi)\rgl|\,,\,
                   |\lgl|\xi|^k,\vphi(\xi)\rgl|\big)\le
      C_1^kC_2C_3\le(C_1C_2C_3)^k,$$
$$\lgl\vphi|\vphi\rgl^k\le(C_2C_3^2)^k,$$
$$\lgl\vphi|\vphi\rgl^{\ga_k}\le(C_2C_3^2)^{k+1}\le
      (C_2^2C_3^4)^{k}.$$
2. The Schwarz inequality yielding
$1=\lgl1,\vphi\rgl\le(\int\limits_L 1)^{\frac{1}{2}}\|\varphi\|_2
=(\mbox{\rm m}(L)\lgl\vphi|\vphi\rgl)^{\frac{1}{2}}$, we conclude
$\lgl\vphi|\vphi\rgl\ge C_2^{-1}$ and from this, in turn,
$$\lgl\vphi|\vphi\rgl^{-k}\le C_2^k,$$
$$\lgl\vphi|\vphi\rgl^{-\ga_k}\le(C_2^2)^k.$$
3. The case $n=0$ being trivial, note that there exists $C_0>1$
such that $k^n\le C_0^k$ for all $k\in\N$. Consequently,
$$|(\bet_k)_n|\le(|\bet_k|+n-1)^n\le(pk+n-1)^n\le(pkn)^n\le(C_0C)^k$$
where $C:=\max(1,(pn)^n)$; the third inequality in the preceding
chain is based on $0\le(pk-1)(n-1)$.
\ep

Now the (ep)-property for
$\eps^{-k}   \rmd^m (P_k^{(3)})_\eps=
               \rmd^m (P_k^{(3)})=\rmd^m(\lgl v_k,\,.\,\rgl)$
is clear from Proposition \rf{listeb} and the remark preceding it.
Regarding
$\eps^{\ga_k}\rmd^m (P_k^{(2)})_\eps=\rmd^m (P_k^{(2)})=
         \rmd^m(\lgl\,.\,|\,.\,\rgl^{\ga_k})$ we
obtain from the chain rule that
$(\rmd^m\lgl\,.\,|\,.\,\rgl^{\ga_k})(\vphi)(\psi_1,\dots,\psi_m)$
is given as a certain sum
of terms of the form
\bea\lb{termdlr2}
   (\ga_k)_l\lgl\vphi|\vphi\rgl^{\ga_k-l}\cdot
   (\rmd^{j_1}\lgl\,.\,|\,.\,\rgl)(\dots)\cdot\dots\ldots\cdot
   \rmd^{j_l}\lgl\,.\,|\,.\,\rgl)(\dots),
\eea
the groups of three dots in parentheses
standing for $\vphi$ and certain subsequences of
$(\psi_1,\dots,\psi_m$). ($j_1+\dots+j_l=m$ and, for any
non-vanishing term of the above form, $j_1,\dots,j_l\in\{1,2\}$.)
(\rf{termdlr2}) immediately allows the application of
Propositions \rf{ebalg} and \rf{listeb}, again in connection with
the remark preceding the latter, thereby establishing the
(eb)-property also for $\eps^{\ga_k}\rmd^m (P_k^{(2)})_\eps$.
For both terms treated so far, the constants occurring in the
(eb)-estimate obviously are independent of $\eps$.
Observe that the case $n=0$ of (\rf{estpk}) is already settled completely
on the basis of the results obtained so far, due to $g$ being
globally bounded on $\R$.

We now turn to the
remaining one of the three terms which, to be sure, is the most
difficult one to handle: We have to show
$(\eps^{n\eta}\rmd^m (P_k^{(1)})_\eps)_k$ to constitute an
(eb)-sequence, with the corresponding constant not depending on
$\eps$. Again according to the chain rule,
the $m$-th differential of
$\eps^{n\eta}g\big(\eps^{-\ga_k}\lgl\vphi|\vphi\rgl^{\ga_k}
e(v(\vphi))\big)$, evaluated at $\vphi;\psi_1,\dots,\psi_m$,
is given as a sum of terms of the form
$$\eps^{n\eta-l\ga_k}g^{(l)}\big(\eps^{-\ga_k}\lgl\vphi|\vphi\rgl^{\ga_k}
e(v(\vphi))\big)\cdot f_1(\dots)\cdot\dots\ldots\cdot
            f_l(\dots)\qquad\qquad(1\le l\le m)$$
(maintaining the convention that a group of three dots at a
differential's argument's place always denotes a certain
subsequence of $(\vphi;\psi_1,\dots,\psi_l)$) where each
$f_{l'}(\dots)$ is of the form
\be\lb{rmdj}
\rmd^j\big(\lgl\,.\,|\,.\,\rgl^{\ga_k}\cdot (e\circ
           v)\big)(\dots)\qquad\qquad(1\le j\le l\le m).
\end{equation}
On the basis of Leibniz' rule, Proposition
\rf{ebalg} and the fact that
$\lgl\vphi|\vphi\rgl^{\ga_k}$ together with all its differentials
is already known to be (eb),
it will suffice to deal with terms of the form 
\be\lb{}
\eps^{n\eta-l\ga_k}g^{(l)}\big(\eps^{-\ga_k}\lgl\vphi|\vphi\rgl^{\ga_k}
   e(v(\vphi))\big)\cdot
   \rmd^{i_1}(e\circ v)(\dots)\cdot\dots\ldots\cdot
   \rmd^{i_l}(e\circ v)(\dots)
\end{equation}
where $0\le i_{l'}\le l$ and $i_1+\dots i_l\le m$.
To this end, we have to analyze $\rmd^{i}(e\circ v)$ for
$0\le i\le j\le l\le m$.
Once more by the chain rule, this is a sum of products of
$e^{(r)}(v(\vphi))$ with $r$ factors which are
differentials of $v$ ($0\le r\le i$). The proof of Proposition
\rf{listeb} shows that $\lgl\vphi|\vphi\rgl^{\frac{1}{2}-r'}$
is bounded on bounded sets for $r'\in\N$. From this it follows that also
the differentials of $v$
are bounded on bounded sets; as they do not depend on $k$,
they form an (eb)-sequence in a trivial manner.
By Proposition \rf{ebalg} again we can discard them for the rest of
the argument.
Thus we are left with estimating
\be\lb{gschiach}
\eps^{n\eta-l\ga_k}g^{(l)}\big(\eps^{-\ga_k}\lgl\vphi|\vphi\rgl^{\ga_k}
e(v(\vphi))\big)\cdot e^{(r_1)}(v(\vphi)) \cdot\dots\cdot
            e^{(r_l)}(v(\vphi))\qquad(1\le l\le m)
\end{equation}
where $0\le r_{t}\le i_{t}\le l$ ($1\le t\le l$)
and $r_1+\dots+r_l\le m$.
Now $e^{(r)}(v(\vphi))$ can be written as
$$e^{(r)}(v(\vphi))=e(v(\vphi))\frac{q_r(v(\vphi))}{v(\vphi)^{2r}}$$
where $q_r$ is a certain polynomial of degree $r-1$. Consequently,
(\rf{gschiach}) takes the form
$$\eps^{n\eta-l\ga_k}g^{(l)}(X)\cdot e(v(\vphi))^l \cdot
            \frac{1}{v(\vphi)^{2n}}\cdot v(\varphi)^{2(n-\bar r)}
            \prod_{t=1}^{l}q_{r_t}(v(\vphi))$$
where we have set $X:=\eps^{-\ga_k}\lgl\vphi|\vphi\rgl^{\ga_k}
e(v(\vphi))$ and $\bar r:=\sum\limits_{t=1}^{l}r_t$, for the sake of
brevity. Now expand $e(v(\vphi))^l$ according to
$$e(v(\vphi))^l=X^{l(1-\frac{\eta}{\ga_k})}\cdot
 (\eps^{\ga_k}\lgl\vphi|\vphi\rgl^{-\ga_k})^{l(1-\frac{\eta}{\ga_k})}
 \cdot e(v(\vphi))^{l\frac{\eta}{\ga_k}}$$
and regroup the terms in the following way as to obtain the desired
estimates:

1.\ Collecting all powers of $\eps$, we obtain
$\eps^{n\eta-l\ga_k}\cdot\left(\eps^{\ga_k}\right)
    ^{l(1-\frac{\eta}{\ga_k})}=\eps^{(n-l)\eta}\le1.$

2.\ For $|X|\le1$, we have
$\left|g^{(l)}(X)\cdot X^{l(1-\frac{\eta}{\ga_k})}\right|\le
   \left|g^{(l)}(X)\right|\le\|g^{(l)}\|_\infty$
(note that $0<\eta\le1<\ga_1=1+1\le\ga_k$ and that,
consequently, $\frac{\eta}{\ga_k}\le\frac{1}{2}$ for all $k\in\N$),
while for $|X|\ge1$ and $c_l$ denoting a positive constant
dominating $|x|^{l+1}|g^{(l)}(x)|$ for all $x\in\R$ (see the
remarks after the introduction of $g$), we obtain
$$\left|g^{(l)}(X)\cdot X^{l(1-\frac{\eta}{\ga_k})}\right|\le
  c_l\left|X\right|^{-l-1+l-l\frac{\eta}{\ga_k}}\le c_l.$$
Altogether, the function
$X\mapsto g^{(l)}(X)\cdot X^{l(1-\frac{\eta}{\ga_k})}$
($l=1,\dots,n$) is globally bounded by a positive constant
larger or equal to 1, say, $C_g$.

3.\ The following term, that is
$\lgl\vphi|\vphi\rgl^{-\ga_k l(1-\frac{\eta}{\ga_k})}=
  \lgl\vphi|\vphi\rgl^{-l(\ga_k-\eta)}$
gives rise to an (eb)-sequence letting $k=1,2,\dots$:
This is immediate from
$\lgl\vphi|\vphi\rgl^{-(\ga_k-\eta)}=
  \lgl\vphi|\vphi\rgl^{-\ga_k}\cdot
  \lgl\vphi|\vphi\rgl^\eta$
and Propositions \rf{listeb} and \rf{ebalg}, together with the
observation that $\lgl\vphi|\vphi\rgl^\eta$ is bounded
on bounded sets. Hence
for a given bounded subset $B$ of $\ca_0(\R)$ there exists
a constant $C_1\ge1$ satisfying
$\lgl\vphi|\vphi\rgl^{-l(\ga_k-\eta)}
   \le C_1^k$
for all $\vphi\in B$, $k\in\N$.

4.\
$e(v(\vphi))^{l\frac{\eta}{\ga_k}}\cdot \frac{1}{v(\varphi)^{2n}}$
can be rewritten as
$$e\left(\frac{\ga_k v(\varphi)}{l\eta}\right)\cdot
   \left(\frac{l\eta}{\ga_k v(\varphi)}\right)^{2n}\cdot
   \frac{1}{(l\eta)^{2n}}\cdot\ga_k^{2n}.$$
Now $e(x)\cdot x^{-2n}$ (with $\frac{0}{0}:=0$) is globally
bounded on $\R$ and $\ga_k^{2n}\le(k+1)^{2n}\le
(2k)^{2n}=4^nk^{2n}$; the latter is (eb) by the proof of
part 3 of Proposition \rf{listeb}. Therefore,
$e(v(\vphi))^{l\frac{\eta}{\ga_k}}\cdot \frac{1}{v(\varphi)^{2n}}$
can be estimated by $C_2^k$ for a suitable constant $C_2\ge1$.

5.\ $\vphi$ ranging over the bounded set $B$ as in
3.\ above, $v(\vphi)$ attains values in a bounded subset of
$\C$. On this set the polynomial
$x^{2(n-\bar r)} \prod\limits_{t=1}^{l}q_{r_t}(x)$
is bounded by some constant $C_3\ge1$.

Summarizing, for any given bounded subset $B$ of $\ca_0(\R)$
there exist constants $C_g,C_1$,\linebreak
$C_2,C_3$ (independent of
$\eps\in I$) such that
$$
\left|\eps^{n\eta-l\ga_k}g^{(l)}\big(\eps^{-\ga_k}\lgl\vphi|\vphi\rgl^{\ga_k}
e(v(\vphi))\big)\cdot e^{(r_1)}(v(\vphi)) \cdot\dots\cdot
            e^{(r_l)}(v(\vphi))\right|\le(C_gC_1C_2C_3)^k
$$
for all $\vphi\in B$. This completes the proof of (\rf{estpk}).
\ep

\subsection{Proof of smoothness of $P$}
Setting $\eps:=1$ in (\rf{estpk}) shows
$(\rmd^n\!P_k)$ to be an (eb)-sequence on
$\ca_0(\R)\times\ca_{00}(\R)^n$, for each $n\in\N_0$.
Theorem \rf{theb} now implies that $P$ as defined in
\rf{defpq} is smooth, that the differentials of
$P$ can be computed term-wise and that all the series for 
$\rmd^n\!P$ ($n\in\N_0$) converge with respect to $\tau_\infty$.
\ep

\subsection{Proof of moderateness of $P$}
Let $B$ be a bounded subset of $\cd(\R)$ and assume
$C_n$ ($n\in\N_0$) to be appropriate constants as to satisfy
(\rf{estpk}) for all $\vphi\in B\cap\ca_0(\R)$,
$\psi_1,\dots,\psi_n\in B\cap\ca_{00}(\R)$.
Choosing $\eta:=1$, say, estimate (\rf{estpk}) results in
$|(\rmd^n(P_k)_\eps)(\vphi)(\psi_1,\dots,\psi_n)|\le
   C_n^k\cdot\eps^{-\frac{1}{k}-n}.$
Multiplying by $\frac{1}{k!}$ and
forming the infinite sum constituting
$\rmd^n\!P_\eps$, we obtain
$$|(\rmd^n\!P_\eps)(\vphi)(\psi_1,\dots,\psi_n)|\le
   (e^{C_n}-1)\cdot\eps^{-1-n},$$
uniformly on $B$ in the sense specified above.
Hence $P$ satisfies the condition 
equivalent to moderateness given in
Theorem 7.12.
\ep

We proceed to prove $P\notin\cn^d$ resp.\ $P\in\cn^e$.
For the former negligibility property, instead of the condition
as given in 7.3, we use the equivalent condition
(once again omitting the quantifiers
``$\forall K\subset\subset\Om\ \forall\al\in\N_0^s$'')
\beas
R\in\cn^d\Leftrightarrow
\forall n\ \exists q\ 
   \forall B\mbox{\rm \,(bounded)\,}\subseteq \ca_q(\R^s):
   \sup\limits_{x\in K,\ \vphi\in B}|\pa^\al(R(S_\eps\vphi,x))|
   =O(\eps^n)
\eeas
occurring as 1$^\circ$ in Theorem 18 of \cite{JEL}.
Observe that for the application of the latter theorem, we need the
fact that $P\in\ce^d_M$ which has been shown above.
For $\cn^e$ we use the modified
defining condition as given previously in this section
(cf.\ the discussion of $\Phi_0$):
\beas
R\in\cn^e\Leftrightarrow
\forall n\ \exists q\ 
   \hphantom{mmmmmm}
   \forall\vphi\in\ca_q(\R^s):\hphantom{m}
   \sup\limits_{x\in K}\hphantom{m}|\pa^\al(R(S_\eps\vphi,x))|
   =O(\eps^n)
\eeas
Clearly, our choice of the above forms of the respective conditions
is motivated by the intention to have them as similar as possible
to highlight the essential difference between them: The estimate on
$|\pa^\al(R^e(S_\eps\vphi,x))|$
is required to hold uniformly on bounded subsets
with respect to $\vphi$ in the former case
as compared to only pointwise in the latter.

\subsection{Proof of $P\not\in\cn^d$}
Set $K:=\{0\}$, $\al:=0$, $n:=1$. We are going to show that for this
set of data the condition for $P$ to belong to $\cn^d$ is violated,
i.e., we are going to show that for every $q\in\N$ there exists a
bounded subset $B$ of $\ca_q(\R)$ such that
$\sup\limits_{\vphi\in B}|(P(S_\eps\vphi,0))|$
is {\it not} of order $O(\eps)$. To this end, let $q\in\N$.
Since $v_\frac{1}{2},v_0,v_1,\dots,v_{q+1}$ are linearly
independent in $\cd'(\R)$ there exist $\vphi_0,\vphi_1\in\ca_q(\R)$
satisfying
\beas \lgl v_\frac{1}{2},\vphi_0\rgl=0,&\qquad&
               \lgl v_{q+1},\vphi_0\rgl=1,\\
      \lgl v_\frac{1}{2},\vphi_1\rgl=1,&\qquad&
               \lgl v_{q+1},\vphi_1\rgl=1.
\eeas
Setting $\vphi_\lambda:= (1-\lambda)\vphi_0+\lambda\vphi_1$
($0\le\lambda\le1$), $B:=\{\vphi_\lambda\mid0\le\lambda\le1\}$
is a bounded subset of $\ca_q(\R)$; moreover,
$\lgl v_\frac{1}{2},\vphi_\lambda\rgl=\lambda$.
For each $\lambda$ in a suitable interval
$(0,\lambda_0]$ we are going to specify
some $\eps_\lambda \in I$ with $\eps_\lambda \to0$ as
$\lambda\to0$ such that
$P_{\eps_\lambda }(\vphi_\lambda,0)\to\infty$ ($\lambda\to0$).
Consequently, $\sup\limits_{\vphi\in B}|(P(S_\eps\vphi,0))|$
is not even of order $O(1)$, i.e. not even bounded as $\eps\to0$.
The (nonnegative) function defined by the assignment
$$\lambda\mapsto\eps_\lambda :=
  \lgl\vphi_\lambda|\vphi_\lambda\rgl
  \cdot e(v(\vphi_\lambda))^{\frac{1}{\ga_{q+1}}}
$$
is continuous for $\lambda\in[0,1]$, strictly positive
for $\lambda>0$ and satisfies $\eps_0=0$.
Hence there exists $\lambda_0>0$
such that $\eps_\lambda\in I$
for $0\le\lambda\le\lambda_0$. Moreover,
$\eps_\lambda \to0$ as $\lambda\to0$.
The general term of the series defining
$P_{\eps}(\vphi_\lambda,0)$ is given by
(apart from the factor $\frac{1}{k!}$)
$$
(P_k)_{\eps}(\vphi_\lambda)=\eps^{-\frac{1}{k}}
   \cdot g\big(\eps^{-\ga_k}\lgl\vphi_\lambda|\vphi_\lambda\rgl^{\ga_k}
           e(v(\vphi_\lambda))\big)\cdot
   \lgl\vphi_\lambda|\vphi_\lambda\rgl^{\ga_k}\cdot
   \lgl v_k,\vphi_\lambda\rgl.
$$
For $k=1,\dots,q$ this expression vanishes identically on $B$ due
to $\vphi_\lambda\in\ca_q(\R)$. For $k\ge q+2$ it can be estimated
by
$\eps^{-\frac{1}{k}}\cdot\frac{1}{2}\cdot C^k$ 
(note that $\|g\|_\infty=\frac{1}{2}$) for some constant $C\ge1$
being independent of $\lambda$ since
$\lgl\vphi|\vphi\rgl^{\ga_k}\lgl v_k,\vphi\rgl$
forms an (eb)-sequence. Consequently,
$$\sum_{k=q+2}^\infty\frac{1}{k!}
  (P_k)_{\eps}(\vphi_\lambda)\le\eps^{-\frac{1}{q+2}}\cdot
  \frac{1}{2}\cdot e^C.$$
It remains to look at the leading term, that is, 
$(P_{q+1})_{\eps}(\vphi_\lambda)$
(again omitting $\frac{1}{(q+1)!}$).
Setting $\eps:=\eps_\lambda $ it takes the value
$$
\eps_\lambda ^{-\frac{1}{q+1}}\cdot
   g(1)\cdot
   \lgl\vphi_\lambda|\vphi_\lambda\rgl^{\ga_{q+1}}\cdot1
   \ =\ \eps_\lambda ^{-\frac{1}{q+1}}\cdot\frac{1}{2}\cdot
   \lgl\vphi_\lambda|\vphi_\lambda\rgl^{\ga_{q+1}}.
$$
Altogether we obtain
\beas
P_{\eps_\lambda }(\vphi_\lambda,0)&\ge&
   \frac{1}{(q+1)!}\cdot
   \eps_\lambda ^{-\frac{1}{q+1}}\cdot\frac{1}{2}\cdot
   \lgl\vphi_\lambda|\vphi_\lambda\rgl^{\ga_{q+1}}
   -\eps^{-\frac{1}{q+2}}\cdot\frac{1}{2}\cdot e^C\\&=&
   \eps_\lambda ^{-\frac{1}{q+1}}\cdot\frac{1}{2}
   \left[\frac{\lgl\vphi_\lambda|\vphi_\lambda\rgl^{\ga_{q+1}}}
   {(q+1)!}-\eps_\lambda ^{\frac{1}{(q+1)(q+2)}}\cdot
   e^C\right]
\eeas
which tends to infinity as $\lambda\to0$ (and, consequently,
$\eps_\lambda \to0$), due to $\lgl\vphi_\lambda|\vphi_\lambda\rgl$
being bounded from below uniformly for $\lambda\in[0,1]$
and the second term in the square
bracket vanishing in the limit.
\ep

\subsection{Proof of $P\in\cn^e$}
Let $K\subset\subset\R$, $\al:=0$ (note that
$P_\eps(\vphi,x)$ does not depend on $x$)
and $n\in\N$ be given; we claim that $q:=n-1$ is an appropriate
choice for showing that
$$\forall\vphi\in\ca_q(\R):
   \sup\limits_{x\in K}|P_\eps(\vphi,x)|
   =O(\eps^n).$$
Let $\vphi\in\ca_q(\R)=\ca_{n-1}(\R)$. If
$\lgl v_\frac{1}{2},\vphi\rgl\le0$ then $v(\vphi)\le0$ and,
consequently, $e(v(\vphi))=0$ which in turn implies
$P_\eps(\vphi,x)$=0 for all $x\in\R$ and all $\eps\in I$.
Thus we may assume
that $\lgl v_\frac{1}{2},\vphi\rgl>0$. But then also
$v(\vphi)$ and, in turn,  $e(v(\vphi))$ are positive.
Taking into account that
$|g(x)|=|\frac{1}{x}\cdot\frac{x^2}{x^2+1}|\le|\frac{1}{x}|$
for $x\ne0$ we obtain the following estimate:
\beas
|g\big(\eps^{-\ga_k}\lgl\vphi|\vphi\rgl^{\ga_k}
   e(v(\vphi)\big)\cdot\eps^{-\ga_k}\lgl\vphi|\vphi\rgl^{\ga_k}\cdot
      \eps^k\lgl v_k,\vphi\rgl|&\le&
\eps^{k}\frac{|\lgl v_k,\vphi\rgl|}{e(v(\vphi))}.
\eeas
Choosing a constant $C$ satisfying $|\lgl v_k,\vphi\rgl|\le C^k$
for all $k\in\N$ (note that $(v_k)_k$ is (eb) by Proposition
\rf{listeb}) we finally arrive at
$$|P_\eps(\vphi,x)|\le\sum_{k=q+1}^{\infty}
  \frac{1}{k!}\cdot\frac{\eps^k C^k}{e(v(\vphi))}\le
  \eps^{q+1}\cdot\frac{e^C}{e(v(\vphi))}
$$
thereby completing the proof of $P\in\cn^e$.
\ep

Now we turn to briefly discussing $Q$. In what follows we will
tacitly assume all (eb)-questions to be handled appropriately.
After scaling $\vphi$ and dropping the factor $\frac{1}{k!}$
the typical term of the series defining $Q$ takes the form
$$\eps^{-\frac{1}{k}}\cdot
   h_k\big(\frac{1}{\eps}\lgl\vphi|\vphi\rgl\,v(\vphi)\big)\cdot
   \lgl\vphi|\vphi\rgl^{\ga_k}\cdot
   \lgl v_k,\vphi\rgl.
$$
As with $P$, $\rmd^m(\lgl\,.\,|\,.\,\rgl^{\ga_k})$
and $\rmd^m(\lgl v_k,\,.\,\rgl)$ are (eb) for all $m\in\N_0$.
Modulo some (eb)-arguments again, the non-trivial part of dealing with
$\rmd^m h_k\big(\frac{1}{\eps}\lgl\vphi|\vphi\rgl\,v(\vphi)\big)$
consists in getting to grips with
$\eps^{-l}h_k^{(l)}\big(\frac{1}{\eps}\lgl\vphi|\vphi\rgl\,v(\vphi)\big)$
for $l\le m$. Thanks to the harmless leading factor $\eps^{-l}$
(as compared to $\eps^{-l\ga_k}$ in the case of $P$)
it is sufficient to note that there exists some constant
$C\ge1$ satisfying $\|h_k^{(l)}\|_\infty\le C^k$
for all $k\in\N$ and $0\le l\le m$
(observe that $\si$ and $g$ are globally bounded
together with all their derivatives). Summarizing,
we obtain that for all $m\le n$ the sequences (with respect to $k$) 
$\eps^m\rmd^m h_k\big(\frac{1}{\eps}\lgl\vphi|\vphi\rgl\,v(\vphi)\big)$
and, consequently,
$$
\eps^{n}\cdot\rmd^n\Big(
   h_k\big(\frac{1}{\eps}\lgl\vphi|\vphi\rgl\,v(\vphi)\big)\cdot
   \lgl\vphi|\vphi\rgl^{\ga_k}\cdot
   \lgl v_k,\vphi\rgl\Big)
$$
are (eb), with the respective constants not
depending on $\eps$. From this, smoothness and moderateness of
$Q$ follow. To obtain the proof of $Q\notin\cn^d$
from the proof of $P\notin\cn^d$ simply replace
the former definition of $\eps_\lambda$ by
$\eps_\lambda :=
  \lgl\vphi_\lambda|\vphi_\lambda\rgl^{\frac{3}{2}}
  \cdot \lambda$
and use the fact that $h_k(1)=\|h_k\|_\infty=1$.
Finally, to show that $Q\in\cn^e$, fix $\vphi\in\ca_q(\R)$.
The case $\lgl v_\frac{1}{2},\vphi\rgl=0$ being trivial,
assume that $\lgl v_\frac{1}{2},\vphi\rgl\ne0$.
For $\eps\le\frac{2}{3}\lgl\vphi|\vphi\rgl^{\frac{3}{2}}
|\lgl v_\frac{1}{2},\vphi\rgl|$ we have
\beas
\quad\ 
\Big|h_k\Big(\frac{1}{\eps}\lgl\vphi|\vphi\rgl^{\frac{3}{2}}
\lgl v_\frac{1}{2},\vphi\rgl\Big)\Big|=
\Big|2g\Big(\frac{1}{\eps}\lgl\vphi|\vphi\rgl^{\frac{3}{2}}
\lgl v_\frac{1}{2},\vphi\rgl\Big)\Big|^{\ga_k}
\le\eps^{\ga_k}\bigg(\frac{2}{\lgl\vphi|\vphi\rgl^{\frac{3}{2}}
\big|\lgl v_\frac{1}{2},\vphi\rgl\big|}\bigg)^{\ga_k}.
\eeas
The rest of the argument is similar to that for $P$.

The reader might ask if it is indeed necessary to come up with
counterexamples as complicated as $P$ and $Q$ certainly are.
The author doubts
that easier ones might be possible. This view is based on reflecting
on the r\^oles each of the three factors constituting
a single term of the series
(for $P$, say) in fact has to play:
\begin{itemize}
\item $\lgl v_k,\vphi\rgl$ distinguishes between the spaces
      $\ca_q(\R)$; this is crucial for the negligibility
      properties.
\item $\lgl\vphi|\vphi\rgl^{\ga_k}=\lgl\vphi|\vphi\rgl^k\cdot
      \lgl\vphi|\vphi\rgl^{\frac{1}{k}}$, on the one hand,
      after scaling of $\vphi$
      compensates for the factor $\eps^k$ generated by scaling
      $\vphi$ in $\lgl v_k,\vphi\rgl$. On the other hand, it
      introduces a factor $\eps^{-\frac{1}{k}}$ making the
      first non-vanishing term of the series the dominant one
      as $\eps\to0$.
\item $g(\lgl\vphi|\vphi\rgl^{\ga_k} e(\lgl v,\vphi\rgl))$
      allows the pointwise vs.\ uniformly distinction being
      necessary to obtain $P\notin\cn^d$, $P\in\cn^e$.
      Though
      $g(\lgl\vphi|\vphi\rgl^{\ga_k}\lgl v,\vphi\rgl)$ would suffice
      to achieve the latter, this alternative choice
      for the argument of $g$ would produce, via the chain rule,
      a factor $\eps^{-n(k+\frac{1}{k})}$ in the $k$-th term of
      $\rmd^n\!P_\eps$ which would be disastrous for the moderateness of
      $P$. The function $e$ (together with $\eps^{-\ga_k}$ in the
      argument of $g$) suppressing this
      unwanted factor, $P$ becomes moderate in the end.
\end{itemize}
Similar arguments apply to $Q$.

\section{Classification of smooth Colombeau algebras between $\cg^d(\Om)$
and $\cg^e(\Om)$}\lb{spec}
\subsection{The development leading from $\cg^e(\Om)$ to $\cg^d(\Om)$}
This section, in fact, does justice to the title of the paper
by going back to the roots of Colombeau algebras constructed
according to the scheme outlined in section 3 of Part I.
Surveying
the range of algebras lying between the algebras
$\cg^d(\Om)$ and (the smooth version of) $\cg^e(\Om)$,
we will discuss, in particular, to which extent at least
the definition of the  algebra $\cg^1(\Om)$
of \cite{CM} (which can be located within
that range) has to be modified to obtain diffeomorphism invariance.
To be sure, the introduction of $\cg^1(\Om)$ has
to be considered as the decisive step towards the 
construction of a diffeomorphism invariant Colombeau algebra.
The result of our analysis will be the construction of a
diffeomorphism invariant Colombeau algebra $\cg^2(\Om)$
which is non-isomorphic to $\cg^d(\Om)$, yet closer to
$\cg^1(\Om)$ than $\cg^d(\Om)$ is.

Apart from $\cg^e(\Om)$, all algebras to be considered in this and
the subsequent section
have $\cc^\infty(U(\Om))$ resp.\
$\cc^\infty(\ca_0(\Om)\times\Om)$ as their basic space. In
particular, they are smooth algebras in the sense that
representatives $R$ have to be smooth also with respect to $\vphi$.
The maps $\si$, $\io$, $D_i$ and the actions induced by a diffeomorphism are
defined as in 7.1 and 5.5--5.8, respectively.
The algebras will differ, however, as to the type of test objects
used for selecting the moderate resp.\ negligible members from the
basic space.
We begin by briefly reviewing the development leading from $\cg^e(\Om)$
via $\cg^1(\Om)$ to $\cg^d(\Om)$.\ms
{\bf Features distinguishing $\cg^1(\Om)$ from $\cg^e(\Om)$:}
\begin{itemize}
\item[(1.0)] Smooth dependence of $R$ on $(\vphi,x)$ rather than
arbitrary dependence on $\vphi$ and smoothness only with respect to
$x$.
\item[(1.1)] Dependence of test objects on $\eps$.
\item[(1.2)] Asymptotically vanishing moments of test objects
as compared to the stronger condition $\phi(\eps)\in\ca_q(\R^s)$
for all $\eps$ (which would be the na\"\i ve analog of
$\vphi\in\ca_q(\R^s)$ in the case of $\cg^e(\Om)$).
\end{itemize}

{\bf Features distinguishing $\cg^d(\Om)$ from $\cg^1(\Om)$:}
\begin{itemize}
\item[(2.1)]
Dependence of test objects also on $x\in\Om$ (in fact, smooth
dependence).
\item[(2.2)] In testing for moderateness, test objects
for $\cg^d(\Om)$ can take
arbitrary values in $\ca_0(\R^s)$, independently of any moment
condition.
\end{itemize}

Let us analyze briefly how compelling the above changes
in the definitions in fact
are if a diffeomorphism invariant algebra is to be obtained.
We refrain from questioning (1.0), i.e.,
smoothness of $R$ with respect to $(\vphi,x)$,
as well as from questioning the smoothness of
test objects with respect to $x$ in the sense of (2.1). Both
properties being used
in the proof of diffeomorphism invariance in an essential way, they
are absolutely necessary from a pragmatic point ot of view to
guarantee the smoothness of
$(\hat\mu R)(S_\eps\ti\phi(\eps,\ti x),\ti x)
     =R(S_\eps\phi(\eps,\mu\ti x),\mu\ti x)$
with respect to $\ti x$ (see section 7), to be
sure. Of course, this does not amount to say that
we have a formal proof that
for the diffeomorphism invariance of an algebra, smoothness
of $R$ with respect to $\vphi$ or of test objects with respect to $x$
are logically necessary.

Smoothness of test objects with respect to $\eps$ definitly is not
an issue of striking importance: The equivalence of conditions
(B) and (C) in Theorem 10.5
(resp.\ of conditions $(\mathrm{B}')$ and $(\mathrm{C}')$
in Theorem 10.6)
shows that test objects of
the form $\phi\in\cc^{[\infty,\Om]}_b(I \times\Om,\ca_0(\R^s))$
give rise to $\cg^d(\Om)$ (via using them for testing
moderateness resp.\ negligibility) independently of the
assumption of smoothness of $\eps\mapsto\phi(\eps,x)$.
With the appropriate respective modifications of the proof,
this statement is valid for all types of test objects
being dependent on $\eps$ or $(\eps,x)$,
that is, it is true for all nine
types $[z,Y]$ where $z$ is one of $\eps x$ or $\eps$ (see below for the
definition of these types).

Next, if for a given diffeomorphism $\mu:\ti\Om\to\Om$
the induced map $\hat\mu:\ce(\Om)\to\ce(\ti\Om)$ is to extend the usual
action $\mu^*$ induced by $\mu$ on distributions then we
necessarily have to set $\hat\mu=\bar\mu^*$ , i.e.,
$(\hat\mu R)(\ti\vphi,\ti x)=R(\bar\mu(\ti\vphi,\ti x))$
with $\bar\mu$ as defined in 5.7.
For purposes of testing, we have, in turn, no other choice than to
consider
$\hat\mu_\eps=\bar\mu_\eps^*$ where
$$
\bar\mu_\eps(\ti\vphi,\ti x)
    =\left(\ti\vphi
    \left(\frac{\mu^{-1}(\eps.+\mu\ti x)-\ti x}{\eps}\right)\cdot
        |\det D\mu^{-1}(\eps.+\mu\ti x)|\,,\,\mu\ti x\right).
$$
From
$(\hat\mu_\eps R)(\ti\vphi,\ti x)=R(\bar\mu_\eps(\ti\vphi,\ti x))$
it is now evident that a moderate (resp.\ negligible)
function $R$ from the basic space
has to accept test objects which are dependent on $\eps$ as well as
on $x$ ((1.1.) and (2.1))
if $\hat\mu R$ is still to be moderate (resp.\ negligible)
(see the discussion preceding Theorem 7.14).
(1.2) is compelling since the property that certain moments
of a test object have to
vanish simply is not invariant under $\bar\mu$ resp.\
$\bar\mu_\eps$. The moments of the transformed test objects only
vanish asymptotically. This has the further consequence that accepting
(2.1) raises the question of how to handle
asymptotically vanishing moments with respect to
uniformity in $x\in\Om$: Since all the definitions and theorems
involve uniformity on compact subsets of $\Om$ it seems reasonable
to adopt this condition also for the asymptotically vanishing moment
property, possibly even for all derivatives
$\pa_x^\al\phi$ of a test object $\phi(\eps,x)$.
We will discuss several variants below.

So there only remains change (2.2) for which there seems to be no
apparent necessity. To be sure,
accepting (2.2) widens the range of permissible test
objects, thereby in turn reducing $\ce_M(\Om)$ and $\cn(\Om)$
in size (see example \rf{R2345}\,(i) below).
Yet it has to be admitted that by this reduction,
no generalized functions which are of interest either in the 
development of the theory or in applications are lost.
Quite to the
contrary, accepting (2.2) has the advantage that the definition of
$\ce_M(\Om)$ becomes simpler and, above all, that considerable
flexibility is gained in how to define $\cn(\Om)$, as respective
glances at
Theorems 7.9 and \rf{thcond0}
reveal. Nevertheless, the preceding
discussion leaves open the possibility that a diffeomorphism invariant
Colombeau algebra $\cg^2(\Om)$ could be constructed
avoiding (2.2). $\cg^2(\Om)$ would be closer to
$\cg^1(\Om)$ than $\cg^d(\Om)$ is;
the preceding considerations seem to suggest that
passing from $\cg^1(\Om)$ to $\cg^2(\Om)$ would represent the 
minimal modification of $\cg^1(\Om)$ leading to a diffeomorphism
invariant Colombeau algebra. In any case, a construction as
envisaged above would yield a second example of a diffeomorphism
invariant Colombeau algebra.

\subsection{Classification of test objects}
The term ``test object'' will always refer to some element of
$\cc^\infty_b(I \times\ca_0(\R^s))$; apart from functions
$\phi(\eps,x)$, this formally also includes test objects of the form
$\phi(\eps)$  (depending only on $\eps$) as well as
elements $\vphi$ of $\ca_0(\R^s)$.
From now on, we will write $\lgl\xi^\al,\vphi(\xi)\rgl$
or even only $\lgl\xi^\al,\vphi\rgl$ in place of
$\int\xi^\al\vphi(\xi)\,d\xi$ for $\vphi\in\cd(\R^s)$,
$\al\in\N_0^s$.
\bd
Let $q\in\N$. A
function $\phi:I \to\cd(\R^s)$ (possibly constant and/or
depending also on other arguments, e.g., on $x\in\Om$) 
is said to have {\rm vanishing moments of order $q$}
if $\lgl\xi^\al,\phi(\eps)(\xi)\rgl=0$ for all $\al\in\N_0^s$
with $1\le|\al|\le q$. It is said to have
{\rm asymptotically vanishing moments of order $q$}
if $\lgl\xi^\al,\phi(\eps)(\xi)\rgl=O(\eps^q)$
for all $\al\in\N_0^s$ with $1\le|\al|\le q$. To which extent
this estimate is assumed to hold uniformly with respect
to, e.g., $x\in\Om$ has to be specified separately
(see below).
\et

A function $\phi$ taking values in $\ca_0(\R^s)$ has
vanishing moments of order $q$ if and only if it takes values in
$\ca_q(\R^s)$, actually.
To obtain a classification of
Colombeau algebras lying in the range
between $\cg^d(\Om)$ and (the smooth version of)
$\cg^e(\Om)$ we introduce the
symbols defined in the following list.
They are meant to refer to test objects or to respective notions of
moderateness and negligibility (based on test objects of the
corresponding type) or, finally, to Colombeau algebras defined
as quotients of the respective spaces of moderate functions.

{\bf Parametrization of test objects:}\nopagebreak\smallskip\\
\renewcommand{\arraystretch}{1.5}
\begin{tabular}{ll}
$[\mathrm{c}]$& test objects being single elements
           (``constants'') of $\ca_0(\R^s)$
           resp.\ $\ca_q(\R^s)$ 
            \\
$[\eps]$& test objects depending only on $\eps\in I$\\
$[\eps x]$& test objects depending on $\eps\in I $ as well as on
                  $x\in\Om$\\[6pt]
\end{tabular}

{\bf Moments of test objects:}\smallskip\\
\begin{tabular}{ll}
$[0]$&        test objects taking values in $\ca_0(\R^s)$
              without any restriction on moments\\
$[\mathrm{A}]$&        test objects having asymptotically vanishing
              moments (this symbol always\\[-8pt]
     &       has to refer to test objects of type
            $[\eps]$)\\
$[\mathrm{V}]$&        test objects having vanishing moments, i.e.,
              taking values
                  in some $\ca_q(\R^s)$
\end{tabular}

\hspace*{7pt}\parbox{14.6cm}{The following symbols make sense only for test objects of
type $[\eps x]$; each of them indicates asymptotically vanishing
moments of test objects, possibly also of their
derivatives $\pa_x^\al\phi(\eps,x)$, with the following respective
specifications:}

\begin{tabular}{ll}
$[\mathrm{A}_\mathrm{l}]$&uniformly on the particular $K\subset\subset\Om$
  on which $R$ is being  tested \\[-.7em]
\hphantom{$[\mathrm{A}_\mathrm{l}]$} & (``locally'')\\
$[\mathrm{A}_\mathrm{g}]$& uniformly on each $L\subset\subset\Om$
        (``globally'')\\
$[\mathrm{A}_\mathrm{l}^\infty]$&all derivatives
        uniformly on the particular $K\subset\subset\Om$
   on which $R$ is being \\[-.7em]
\hphantom{$[\mathrm{A}_\mathrm{l}]$} & tested\\
$[\mathrm{A}_\mathrm{g}^\infty]$&all derivatives
        uniformly on each $L\subset\subset\Om$
\end{tabular}

If the compact set $K$ on which $R$ is being tested and/or the
order $q$ of the (asymptotic) vanishing of moments is to be
specified, $K$ resp.\ $q$ will be put as subscript(s)
to the corresponding
A-symbol, e.g., $[\mathrm{A}_\mathrm{l}]_{K,q}$.
Parametrization symbols may be combined with (suitable)
moment symbols. If in a composed symbol $[z,Y]$ $Y$ is one of the A-symbols
then $z=\eps$ resp.\ $z=\eps x$, being redundant, will be omitted
frequently.

Obviously, $[\mathrm{A}_\mathrm{g}^\infty]_q$ implies
$[\mathrm{A}_\mathrm{l}^\infty]_{K,q}$ (for any
$K\subset\subset\Om$)
and $[\mathrm{A}_\mathrm{g}]_q$;
each of the latter, in turn, implies
$[\mathrm{A}_\mathrm{l}]_{K,q}$.
As the examples below show, none of the
reverse implications is true.
\bexs\lb{extestobj}
\begin{itemize}
\item[(i)] Let $\Om:=\R$ and $K:=[-1,+1]$. Define
$\phi_1(\eps,x):=\vphi+\eps^q\sin(x|\ln\eps|)\psi$
where $\vphi\in\ca_q(\R)$ and $\psi\in\ca_{00}(\R)$ with
$\lgl\xi^k,\psi(\xi)\rgl=\de_{kq}$ for $k=1,\dots,q$.
$\lgl\xi^q,\pa_x\phi_1(\eps,0)(\xi)\rgl=\eps^q|\ln\eps|$
is not of order $\eps^q$, yet every
$\pa_x^n\phi_1(\eps,x)$ has bounded image.
Hence $\phi_1$ is of type
$[\mathrm{A}_\mathrm{g}]_q$ and of type
$[\mathrm{A}_\mathrm{l}]_{K,q}$, yet neither
of type $[\mathrm{A}_\mathrm{g}^\infty]_q$ nor of type
$[\mathrm{A}_\mathrm{l}^\infty]_{K,q}$.
\item[(ii)] Let $K\subset\subset\Om$ and set
$\phi_2(\eps,x):=\lambda(x)\vphi_1+(1-\lambda(x)\vphi_2)$
where $\vphi_1\in\ca_q(\R^s)$,
$\vphi_2\in \ca_0(\R^s)\setminus\ca_1(\R^s)$ and
$\lambda\in\cd(\Om)$ such that $0\le\lambda\le1$
and $\lambda\equiv1$ on an open neighborhood of $K$.
Then for $q\in\N$,
$\phi_2$ is both of types
$[\mathrm{A}_\mathrm{l}]_{K,q}$ and
$[\mathrm{A}_\mathrm{l}^\infty]_{K,q}$, yet
neither of type $[\mathrm{A}_\mathrm{g}]_q$
nor of type $[\mathrm{A}_\mathrm{g}^\infty]_q$.
\end{itemize}
\et
\subsection{Classification of full smooth Colombeau algebras}
In the following, we will use the symbols introduced above to
classify smooth Colombeau algebras with respect to the type of test objects
used for testing moderateness.
From a
combinatorial point of view, there are eleven
ways of performing this test, each corresponding to one of the
eleven types of test objects.
The following diagram displays these variants and the
relations between them.
The arrows 
are to be read as implications between the corresponding
notions of moderateness or as inclusion relations between the
corresponding sets of moderate functions (and similarly, for
negligibility, as far as types $[\mathrm{A}]$ and $[\mathrm{V}]$
are concerned). They are {\it not}
representing implications between the properties of test
objects being of the particular types; a diagram of the latter
kind would have to have the arrows reversed, of course.

$$\begin{array}{ccccccc}
{}[\eps x,0]&  &\to&   &[\eps,0]&\to&[\mathrm{c},0]\\
\downarrow&&&&&&\\
{}[\eps x,\mathrm{A}_\mathrm{l}]&\to&[\eps x,\mathrm{A}_{\mathrm{g}}]&&\downarrow&&\\
\downarrow&&\downarrow&&&&\downarrow\\
{}[\eps x,\mathrm{A}_\mathrm{l}^\infty]&\to&[\eps x,\mathrm{A}^\infty_{\mathrm{g}}]
   &\to&[\eps,\mathrm{A}]&&\\
&&\downarrow&&\downarrow&&\\
&&[\eps x,\mathrm{V}]&\to&[\eps,\mathrm{V}]&\to&[\mathrm{c,V}]
\end{array}$$

From the diagram, a useful extension of the characterizations of
$R\in\cn^d(\Om)$
obtained so
far\footnote{$(3^\circ)$, $(4^\infty)$ in 7.9;
$(0^\circ)$--$(2^\circ)$ in \rf{thcond0};
$(\mathrm{A}')$--$(\mathrm{Z}')$ in 10.6;
$(\mathrm{C}'')$, $(\mathrm{Z}'')$ in 10.7 
(in each case assuming $R\in\ce_M^d(\Om)$, in addition).}
is immediate: Test objects in
condition $(4^\infty)$
are of type $[\mathrm{A}_\mathrm{l}^\infty]_{K}$
where $K$ is the compact set
on which $R$ is tested.
This dependence on $K$
of the class of admissible test
objects (going back to \cite{JEL}, Theorem 18)
might seem undesirable since this
class rather ought to be defined universally for tests on
arbitrary compact sets.
However,
it is clear from the diagram that
our list of equivalent conditions
could be extended by adding a further
condition $(5^\infty)$,
obtained from $(4^\infty)$ by replacing
$[\mathrm{A}_\mathrm{l}^\infty]$
with
$[\mathrm{A}^\infty_{\mathrm{g}}]$: 
It suffices to observe that $(3^\circ)$ and $(4^\infty)$
are based on test objects of types
$[\eps x,\mathrm{V}]$ and $[\eps x,\mathrm{A}_\mathrm{l}^\infty]$,
respectively.
It will be a consequence of Corollary \rf{ableitungengratis} below
that a further extension by an analogous equivalent condition $(6^\circ)$,
referring to type $[\mathrm{A}_{\mathrm{g}}]$, can be achieved.

One more glance at the diagram allows to clearly identify the obstacle
against the diffeomorphism invariance of the algebra $\cg^1(\Om)$
defined in \cite{CM}:
The Lemma in section 3 of that article only shows
the $\mu$-transform of test objects of type
$[\eps,\mathrm{A}]$ (being
used in defining $\cg^1(\Om)$)
to be of type $[\eps x,\mathrm{A}_{\mathrm{g}}]$; yet, this is
not sufficient for a positive outcome of the $\mu$-transform of
$R$ being tested for, say,
moderateness, provided $R$ is assumed to be moderate.
So it is Theorem {\bf (T6)} of the blueprint outlined in section
3 which fails for $\cg^1(\Om)$.

Now, if $[X]$ and $[Y]$ are chosen from the set of the eleven types
such that $[Y]$ is located ``south to east'' with respect to $[X]$
in the diagram above (i.e., if $\ce_M[X]\subseteq\ce_M[Y]$) and if,
in addition, $[Y]$ is one of the types $[\mathrm{A}]$ or
$[\mathrm{V}]$ then it easily checked that $\ce_M[X]$ is an algebra
({\bf (T2)}) containing $\cn[Y]\cap\ce_M[X]$ as an ideal ({\bf
(T3)}). Consequently, $\ce_M[X]\big/(\cn[Y]\cap\ce_M[X])$ is an
algebra. We shall refer to algebras arising in this way by the term
``Colombeau-type algebras''. Altogether there are 46 admissible
choices of pairs $[X],[Y]$. In the following definition, we will
specify eleven algebras of this kind, one for each type of
moderateness. These will be the only ones we are to deal with in
the sequel. They might be called ``primary'' since each of the
remaining Colombeau-type algebras can be obtained as some subalgebra
or some quotient algebra of one of them. Note, however, that the
collection of these eleven algebras is not minimal in this respect
(see Theorem 7.10).
\bd\lb{colalg}
If $[X]$ is one of the types $[\mathrm{V}]$ or $[\mathrm{A}]$
define
$$\cg[X]:=\ce_M[X]\big/\cn[X];$$
for types $[0]$ define
\beas
  \cg[\eps x,0]&:=&\ce_M[\eps x,0]\big/
    \big(\cn[\eps x,\mathrm{A}_\mathrm{l}^\infty]\cap\ce_M[\eps
    x,0]\big),\\
  \cg[\eps,0]&:=&\ce_M[\eps,0]\big/\big(\cn[\eps,\mathrm{A}]
    \cap\ce_M[\eps,0]\big),\\
  \cg[\mathrm{c},0]&:=&\ce_M[\mathrm{c},0]\big/\big(\cn[\mathrm{c},
  \mathrm{V}]\cap\ce_M[\mathrm{c},0]\big).
\eeas
(The open set $\Om$ has been omitted from the notation of the respective
algebras.)
\et
We will refer to $\cg[X]$ also by ``the algebra of type $[X]$''.
Each of the algebras mentioned at the beginning of this section
is one of the eleven algebras just defined:
Denoting by $\cg^e_0(\Om)$ the ``smooth part'' of $\cg^e(\Om)$,
i.e., the subalgebra formed by all members having a smooth representative
$R\in\cc^\infty(U^e(\Om))$, it is easy to see that
$\cg^e_0(\Om)=\cg[\mathrm{c},\mathrm{V}]$.
$\cg^1(\Om)$ obviously is equal to $\cg[\eps,\mathrm{A}]$; the
algebra $\cg^2(\Om)$ to be introduced in the following section
is obtained as $\cg[\eps x,\mathrm{A}_\mathrm{g}^\infty]$.
$\cg^d(\Om)$, finally, is given as $\cg[\eps x,0]$.
Observe that according to
Theorem 7.9, $\cn[\eps
x,\mathrm{A}_\mathrm{l}^\infty]$ can be
replaced by $\cn[\eps x,\mathrm{V}]$ in the definition of
$\cg[\eps x,0]$. Moreover, it should be clear from
Examples 7.7 and the discussion preceding them why $\cg[\eps x,0]$
has {\it not} been defined as the quotient 
with respect to
$\cn[\eps x,\mathrm{A}_\mathrm{l}]\cap\ce_M[\eps x,0]$:
This choice (corresponding to using condition $(4^\circ)$ of
\cite{JEL}, Theorem 18)
would invalidate part (iii) of {\bf (T1)}
and thus prevent $\io$ from preserving the product of
smooth functions.

Corollary \rf{ableitungengratis} below will show that
test objects of types $[\mathrm{A}_{\mathrm{g}}]$ and
$[\mathrm{A}^\infty_{\mathrm{g}}]$, respectively,
give rise to the same moderate resp.\ negligible functions.
Moreover, it will follow from Theorem \rf{JT17A} that also
test objects of type $[\mathrm{A}_\mathrm{l}^\infty]$ lead to the same
respective notions of moderateness and negligibility
as test objects of type
$[\mathrm{A}^\infty_{\mathrm{g}}]$ do.
This actually leaves us with nine possibly different
algebras. 

The diagram formed by
the canonical homomorphisms between the resulting nine algebras
is not isomorphic to the previous diagram: On the one hand, as mentioned
above,
$[\mathrm{A}_\mathrm{g}]$, $[\mathrm{A}_\mathrm{l}^\infty]$ and
$[\mathrm{A}_\mathrm{g}^\infty]$ have to be merged
to represent $\cg^2(\Om)=\cg[\eps x,\mathrm{A}_\mathrm{g}^\infty]$.
On the other hand, there is no canonical
homomorphism from $\cg^d(\Om)=\cg[\eps x,0]$ into
$\cg[\eps x,\mathrm{A}_\mathrm{l}]$
since
$\cn[\eps x,\mathrm{A}_\mathrm{l}]\cap\ce_M[\eps x,0]$---not
containing any of $R(\vphi,x):=\lgl\xi^\bet,\vphi(\xi)\rgl$---is
strictly smaller than
$\cn[\eps x,\mathrm{A}_\mathrm{l}^\infty]\cap\ce_M[\eps x,0]$.
We do have canonical homomorphisms, however, both from
$\cg^d(\Om)=\cg[\eps x,0]$ and from $\cg[\eps x,\mathrm{A}_\mathrm{l}]$
into $\cg^2(\Om)=\cg[\eps x,\mathrm{A}_\mathrm{g}^\infty]$.
So we finally arrive at
$$\begin{array}{ccccccc}
{}&  &\cg[\eps x,0]&\to&\cg[\eps,0]&\to&\cg[\mathrm{c},0]\\
&&\downarrow&&\downarrow&&\\
{}\cg[\eps x,\mathrm{A}_\mathrm{l}]&\to&\cg[\eps x,\mathrm{A}^\infty_{\mathrm{g}}]
   &\to&\cg[\eps,\mathrm{A}]&&\downarrow\\
&&\downarrow&&\downarrow&&\\
&&\cg[\eps x,\mathrm{V}]&\to&\cg[\eps,\mathrm{V}]&\to&\cg[\mathrm{c,V}]
\end{array}$$

When establishing {\bf (T1)}--{\bf (T8)} for $\cg^2(\Om)$
in the following section
we will survey briefly which of these theorems is true
for each of the (seven) algebras apart from
$\cg^2(\Om)=\cg[\eps x,\mathrm{A}_\mathrm{g}^\infty]$ and
$\cg^d(\Om)=\cg[\eps x,0]$. Let us anticipate at this point
the facts concerning {\bf (T7)} and {\bf (T8)}, i.e.,
diffeomorphism invariance:                                     
The following counterexamples of moderate functions $R$
for which $\hat\mu R$ fails to be moderate for some diffeomorphism
$\mu$ definitely eliminate six of the
nine algebras from the class of
possibly diffeomorphism invariant ones, beyond any pragmatic reasoning
regarding techniques of proof.
\bexs\lb{R01} Let $\Om:=\R$.
\begin{itemize}
\item[(i)]
The example
$R_0(\vphi,x):=\exp(i\exp(\lgl\vphi|\vphi\rgl))$ presented
in \cite{JEL} shows that all three algebras of type
$[\eps,Y]$ ($Y$=0,A,V), as well as the one of type $[\mathrm{c},0]$
are {\it not} diffeomorphism invariant.
\item[(ii)]
Define
$R_1(\vphi,x):=\lgl\xi,\vphi(\xi)\rgl\cdot\exp(\lgl\vphi|\vphi\rgl)$.
Since $R_1$ vanishes on $\ca_1(\R)\times\R$, it is moderate
with respect to any type $[z,\mathrm{V}]$.
Under the action induced by the diffeomorphism
$\mu(x):=x+e^x$ of $\R$ onto itself,
$R_1$ is transformed to a function $\hat\mu R_1$ which is not moderate with
respect to any type $[z,\mathrm{V}]$ since the values attained by
$$(\hat\mu R_1)(S_\eps\vphi,x)=
  \exp\Big(\frac{1}{\eps}\int\frac{|\vphi(\xi)|^2}{1+e^xe^{\eps\xi}}
  \,d\xi\Big)\cdot
  \Big[\eps\int\xi\vphi(\xi)\,d\xi+
  e^x\int(e^{\eps\xi}-1)\vphi(\xi)\,d\xi\Big]$$
are not of any order $\eps^{-N}$ ($n\in\N$)
even for simple test objects of the
form $\vphi\in\ca_N(\R)$. Therefore, $\hat\mu R_1$
does not pass the test for moderateness.
This example excludes all types
$[z,\mathrm{V}]$ from the class of diffeomorphism
invariant algebras.
\end{itemize}
The details are left to the reader.
\et

Thus we are left with the algebras of types $[\eps x,0]$,
$[\eps x,\mathrm{A}_\mathrm{g}^\infty]$ (together with the two
equivalent types mentioned above) and
$[\eps x,\mathrm{A}_\mathrm{l}]$ as possible candidates for being
diffeomorphism invariant. $[\eps x,0]$ giving rise to the algebra
$\cg^d(\Om)$ introduced in section 7, we will
define $\cg^2(\Om)$ in the following section on the basis of type
$[\eps x,\mathrm{A}_\mathrm{g}^\infty]$ and
prove it to be a diffeomorphism invariant
Colombeau algebra by establishing the corresponding
Theorems {\bf (T1)}--{\bf (T8)}.

For a discussion of
$\cg[\mathrm{A}_\mathrm{l}]$, finally,
we refer to the following section. It is clear from
Example 7.7 that this algebra cannot
be counted among the class of Colombeau algebras due to 
its multiplication not reproducing the product of smooth functions;
moreover, we have to leave it open if it is a differential algebra
at all since we do not know if $\cn[\mathrm{A}_\mathrm{l}]$
is invariant under differentiation.
Nevertheless,
the spaces of moderate resp.\ negligible functions
obtained from type $[\eps x,\mathrm{A}_\mathrm{l}]$ 
test objects turn out to be diffeomorphism
invariant. 
Despite the obvious faults of $\cg[\mathrm{A}_\mathrm{l}]$,
we have included type
$[\mathrm{A}_\mathrm{l}]$ in our scheme,
mainly to allow for a thorough
discussion of condition $(4^\circ)$ of \cite{JEL}, Theorem 18.

Summarizing, the results of this and the
following section show that
$\cg^d(\Om)$ and $\cg^2(\Om)$ are the only diffeomorphism
invariant Colombeau algebras among the eleven (resp.\ nine)
algebras defined in \rf{colalg}.

To complete this section, it remains
to prove that the tests based on
types $[\mathrm{A}^\infty_{\mathrm{g}}]$
and $[\mathrm{A}_{\mathrm{g}}]$ are in fact equivalent.
As demonstrated by Example \rf{extestobj}\,(i),
there are test objects of type $[\mathrm{A}_{\mathrm{g}}]_q$
failing to be of type $[\mathrm{A}^\infty_{\mathrm{g}}]_q$.
Nevertheless, both these classes of test objects do
give rise to the same moderate resp.\ negligible
functions.
This fact will
emerge as an immediate corollary from the following theorem.

\bt\lb{typetype}
Let $\phi\in\cc^\infty_b(I \times\Om,\ca_0(\R^s))$ and let $2\le
q\in\N$. If $\phi$ is of type $[\mathrm{A}_{\mathrm{g}}]_q$ then it also is
of type $[\mathrm{A}^\infty_{\mathrm{g}}]_{q-1}$.
\et
For the proof we need two lemmas.
\blem\lb{type1}
Let $c:I \times\Om\to\R$ have second partial derivatives
$\pa^2_ic$ for some $i\in\{1,\dots,s\}$ ($\pa_i=\frac{\pa}{\pa
x_i}$). Let $q>0$, $0\le r<1$
and assume that $K\subset\subset L\subset\subset\Om$.
If $\sup\limits_{x\in L}|c(\eps,x)|=O(\eps^q)$ and
$\sup\limits_{x\in L}|\pa_i^2c(\eps,x)|=O(\eps^{rq})$
then $\sup\limits_{x\in
K}|\pa_ic(\eps,x)|=O(\eps^{\frac{1+r}{2}q})$.
\et
\pr We consider values of $\eps\in I $ which are less than
$\dist(K,\pa L)$; set $p:=q\frac{1-r}{2}$. For $x\in K$, 
$x+\eps^pe_i\in L$. Taylor's Theorem yields
\beas
c(\eps,x+\eps^pe_i)
   &=&c(\eps,x)+\eps^p\pa_ic(\eps,x)+
      \eps^{2p}\frac{1}{2}\pa_i^2c(\eps,x_\theta)
\eeas
where $x_\theta=x+\theta\eps^pe_i$ for some $\theta\in(0,1)$;
note that also $x_\theta\in L$.
Consequently,
\beas
&\pa_ic(\eps,x)
   =\eps^{-p}
       \underbrace
         {\left(c(\eps,x+\eps^pe_i)
          -c(\eps,x)\right)}
         _{O(\eps^q)}
       -\eps^p\underbrace
         {\frac{1}{2}\pa_i^2c(\eps,x_\theta)}
         _{O(\eps^{rq})}
=O(\eps^{\frac{1+r}{2}q}), &
       \eeas
uniformly for $x\in K$.
\ep\par For the second lemma, we inductively define a sequence of
numbers $r_k$  by setting $r_1:=0$, $r_{k+1}:=\frac{(1+r_k)^2}{4}$
($k\in\N$). Being strictly increasing and bounded by $1$, this
sequence is convergent, its limit being equal to $1$.
\blem\lb{type2}
For every $k\in\N$ the following holds:
Let $c:I \times\Om\to\R$ be smooth with respect to the variable
$x_i$ ($x=(x_1,\dots,x_s)\in\Om$)
for some $i\in\{1,\dots,s\}$. Let $q>0$
and $K\subset\subset L\subset\subset\Om$.
If $\sup\limits_{x\in L}|c(\eps,x)|=O(\eps^q)$ and
$\sup\limits_{x\in L}|\pa_i^mc(\eps,x)|=O(1)$ for all $m\in\N$
then
$\sup\limits_{x\in K}|\pa_ic(\eps,x)|=O(\eps^{\frac{1+r_k}{2}q})$.
\et
\pr
Proceeding by induction, the case $k=1$ is immediate from
Lemma \rf{type1} by setting $r:=r_1=0$. Assume the statement of
the lemma to be true for a particular $k\in\N$. Let $c,i,q,K,L$ be
as specified. Choose $K_1,K_2$ as to satisfy
$K\subset\subset K_1\subset\subset K_2\subset\subset L$.
From $\sup\limits_{x\in L}|c(\eps,x)|=O(\eps^q)$ and
$\sup\limits_{x\in L}|\pa_i^mc(\eps,x)|=O(1)$ for all $m\in\N$
we deduce, by assumption,
$\sup\limits_{x\in K_2}|\pa_ic(\eps,x)|=O(\eps^{\frac{1+r_k}{2}q})$.
Applying the statement of the lemma (for the particular value of
$k$ under consideration) once more, this time
to the function $\pa_ic$, with $\frac{1+r_k}{2}q$ in place of $q$
and for the pair $K_1,K_2$ of compact sets,
we obtain
$\sup\limits_{x\in K_1}|\pa_i^2c(\eps,x)|=O(\eps^{(\frac{1+r_k}{2})^2q})$.
In a last step, we apply Lemma \rf{type1} to conclude
that
$\sup\limits_{x\in K}|\pa_ic(\eps,x)|=O(\eps^{\bar rq})$
where $\bar r=\frac{1}{2}(1+\frac{(1+r_k)^2}{4})=\frac{1+r_{k+1}}{2}$,
thereby showing the statement of the lemma to be true also for
$k+1$.
\ep

{\bf Proof of Theorem \rf{typetype}. }
Let $\phi\in\cc^\infty_b(I \times\Om,\ca_0(\R^s))$ be of type
$[\mathrm{A}_{\mathrm{g}}]_q$ where $2\le q\in\N$.
Denoting $\lgl\xi^\al,\phi(\eps,x)(\xi)\rgl$ by
$c_\al(\eps,x)$ ($\al\in\N_0^s$), we have to show that
$$\sup\limits_{x\in K}|\lgl\xi^\al,\pa_x^\bet\phi(\eps,x)(\xi)\rgl|=
\sup\limits_{x\in K}|\pa^\bet c_\al(\eps,x)|=O(\eps^{q-1})
$$
for $1\le|\al|\le q-1$ and all $K\subset\subset\Om$,
$\bet\in\N_0^s$. Fix $\al\in\N_0^s$ satisfying
$1\le|\al|\le q$.
By assumption, we have
$\sup\limits_{x\in L}|c_\al(\eps,x)|=O(\eps^q)$ 
and $\sup\limits_{x\in L}|\pa^\bet c_\al(\eps,x)|=O(1)$
for all $L\subset\subset\Om$ and all $\bet\in\N_0^s$.
Since $q\frac{1+r_k}{2}\to q$ as $k\to\infty$, Lemma
\rf{type2} yields that
$\sup\limits_{x\in K}|\pa_ic_\al(\eps,x)|=O(\eps^{q-\frac{1}{2}})$
for every $K\subset\subset\Om$ and any $i=1,\dots,s$.
Noting that also
$(q-\frac{1}{2})\frac{1+r_k}{2}\to(q-\frac{1}{2})$,
the same argument, applied to $\pa_i c_\al$ and $\pa_j$
($j=1,\dots,s$)
in place of $c_\al$ and $\pa_i$, respectively, shows that
$\sup\limits_{x\in K}|\pa_j\pa_ic_\al(\eps,x)|=
   O(\eps^{q-(\frac{1}{2}+\frac{1}{4})})$, again
for every $K\subset\subset\Om$ and any $i,j=1,\dots,s$.
By induction, we obtain
$\sup\limits_{x\in K}|\pa^\bet c_\al(\eps,x)|=
   O(\eps^{q-q_\bet})$ for all $\bet\in\N_0^s$
where $q_\bet=\sum\limits_{i=1}^{|\bet|}2^{-i}<1$.
From this we finally conclude that
$\sup\limits_{x\in K}|\pa^\bet c_\al(\eps,x)|=
   O(\eps^{q-1})$
   for all $\bet\in\N_0^s$ and all
$K\subset\subset\Om$.
\ep

\bc\lb{ableitungengratis}
Let $R\in\ce(\Om)$. $R$ is moderate (resp.\ negligible)
with respect to type $[\mathrm{A}_{\mathrm{g}}]$ if and only if it is
moderate (resp.\ negligible)
with respect to type $[\mathrm{A}^\infty_{\mathrm{g}}]$.
\et
\pr
Necessity of the condition being obvious, let us show
sufficiency.
Assuming $R$ to be moderate with respect to type
$[\mathrm{A}^\infty_{\mathrm{g}}]$, fix $\al\in\N_0^s$, $K\subset\subset\Om$.
Choose $N_1\in\N$ such that
$\pa^\al(R(S_\eps\phi_1(\eps,x),x))=O(\eps^{-N_1})$ holds
for every test object $\phi_1$ of type
$[\mathrm{A}^\infty_{\mathrm{g}}]_{N_1}$, uniformly on $K$. Now set $N:=N_1+1$
and pick a test object $\phi$ of type $[\mathrm{A}_{\mathrm{g}}]_N$.
By Theorem \rf{typetype}, $\phi$ is of type
$[\mathrm{A}^\infty_{\mathrm{g}}]_{N-1}$, i.e., of type
$[\mathrm{A}^\infty_{\mathrm{g}}]_{N_1}$. Due to our choice of $N_1$,
$\pa^\al(R(S_\eps\phi(\eps,x),x))=O(\eps^{-N_1})$ resp.\
$O(\eps^{-N})$ follow.
A similar argument applies to negligibility of $R$.
\ep

\section{The algebra $\cg^2$; classification results}\lb{g2}
The algebra $\cg^2(\Om)$ 
of type $[\eps x,\mathrm{A}_\mathrm{g}^\infty]$
to be analyzed below
results from the algebra $\cg^1(\Om)=\cg[\eps,\mathrm{A}]$
of \cite{CM} by applying the minimal modification necessary to
obtain diffeomorphism invariance.
Recall that a test object
$\phi\in\cc_b^\infty(I \times\Om,\ca_0(\R^s))$
is said to be of type $[\eps x,\mathrm{A}_{\mathrm{g}}^\infty]_q$
if $\sup_{x\in K}|\lgl\xi^\al,\pa_x^\bet\phi(\eps,x)(\xi)\rgl|
=O(\eps^q)$ for every $K\subset\subset\Om$, $\bet\in\N_0^s$ and
$\al\in\N_0^s$ with $1\le|\al|\le q$.
Moderateness resp.\ negligibility of $R\in\ce(\Om)=\cc^\infty(U(\Om))$ 
are defined as follows (where $K\subset\subset\Om$ and $\al\in\N_0^s$):
\bd
$R\in\ce(\Om)$ is moderate with respect to type
$[\eps x,\mathrm{A}_\mathrm{g}^\infty]$ if the following condition is
satisfied:
$$\ba{l}
 \forall K\ \forall\al\ \exists N\in\N\
   \forall\phi\in\cc_b^\infty(I \times\Om,\ca_0(\R^s))
   \ \mbox{\rm {which are of type}}\ 
   [\eps x,\mathrm{A}_{\mathrm{g}}^\infty]_N:\\[5pt]
   \hphantom{\forall K\subset\subset\Om\ \forall\al\in\N_0^d\ \exists N}
   \sup\limits_{x\in K}|\pa^\al(R(S_\eps\phi(\eps,x),x))|
   =O(\eps^{-N}).
\ea$$
\et
\bd
$R\in\ce(\Om)$ is negligible with respect to type
$[\eps x,\mathrm{A}_\mathrm{g}^\infty]$ if the following condition is
satisfied:
$$\ba{l}
 \forall K\ \forall\al\ \forall n\in\N\ \exists q\in\N\
   \forall\phi\in\cc_b^\infty(I \times\Om,\ca_0(\R^s))
   \ \mbox{\rm {which are of type}}\ 
   [\eps x,\mathrm{A}_{\mathrm{g}}^\infty]_q:\\[5pt]
   \hphantom{\forall K\subset\subset\Om\ \forall\al\in\N_0^d\ \exists N}
   \sup\limits_{x\in K}|\pa^\al(R(S_\eps\phi(\eps,x),x))|
   =O(\eps^n).
\ea$$
\et
Since we are dealing with $\cg^2(\Om)$ exclusively in the following,
we simply denote the sets of moderate resp.\
negligible functions in the sense of the preceding definitions
by $\ce_M(\Om)$, $\cn(\Om)$.
To establish $\cg^2(\Om)$ as a diffeomorphism invariant
Colombeau algebra we have to convince ourselves that Theorems
{\bf (T1)}--{\bf (T8)} of the scheme presented in section
3 are true on the basis of the preceding
definitions (compare section 7 for the
detailed elaboration of these theorems in the case of
$\cg^d(\Om)$). Though our main interest will be focused on type
$[\mathrm{A}_\mathrm{g}^\infty]$, of course, for each of
{\bf (T1)}--{\bf (T8)} we will specify for which of the
remaining types (apart from $[\eps x,\mathrm{A}_\mathrm{g}^\infty]$
and $[\eps x,0]$)
it holds as well.

To start with, (i) and (ii) of {\bf (T1)}
follow from the corresponding statements with respect to
$\cg^d(\Om)$ (7.4 ,(i),(ii)) for all types
since $[\eps x,0]$ generates the smallest one of all
spaces $\ce_M[X]$. We already know from Example
7.7 that (iii) of {\bf (T1)} is not satisfied for
type $[\mathrm{A}_\mathrm{l}]$. For all the remaining types,
however, the corresponding statement follows immediately 
from part (iii) of 7.4 by
observing that
$\cn[\eps x,\mathrm{V}]\cap\ce_M[\eps x,0]=
\cn[\eps x,\mathrm{A}_\mathrm{l}^\infty]\cap\ce_M[\eps x,0]=
\cn[\eps x,\mathrm{A}_\mathrm{g}^\infty]\cap\ce_M[\eps x,0]$
is contained in each space $\cn[Y]$ where $[Y]$ is different
from $[\mathrm{A}_\mathrm{l}]$.
(The preceding equalities are due to Theorem 7.9
resp.\ to $(4^\infty)\!\Leftrightarrow\!(5^\infty)$ derived in the
preceding section.)
Finally, the proof of part (iv) of {\bf (T1)} given in section
7 for $\cg^d(\Om)$ uses test objects of type
$[\mathrm{c},\mathrm{V}]$ (generating the largest one of all spaces
$\cn[X]$) and therefore is valid for all types.

Theorems {\bf (T2)} and {\bf (T3)} are immediate from
Leibniz' rule for all types.

As it had been the case for $\cg^d(\Om)$, {\bf (T4)}--{\bf (T6)} 
are the hard ones to prove also for $\cg^2(\Om)$.
Fortunately, {\bf (T6)} can be taken
from section 7 with only a slight modification,
as we will see.
For {\bf (T4)} and {\bf (T5)}, however, we need analogs
of Theorems 7.12 and 7.13
for type $[\mathrm{A}_\mathrm{g}^\infty]$ allowing to
express moderateness resp.\ neglibility of $R$ in terms of differentials
of $R_\eps$. To this end, we have to introduce
appropriate classes of sets corresponding to the bounded subsets
$B\subseteq\cd(\R^s)$ occurring in Theorems
7.12 and 7.13.
%
For any closed affine subspace $E_1$ of a locally convex space $E$, let
$\cc^\infty_b(I ,E_1)$ denote the set of all smooth maps
$\vphi:I \to E_1$ having bounded image.
\bd
Let $k\in\N_0$, $q\in\N$.\\ A {\mbox\rm $(k,q)$-class} is a subset
$\cb$ of $\cc^\infty_b(I ,\ca_0(\R^s))\times
\big[\cc^\infty_b(I ,\ca_{00}(\R^s))\big]^k$ satisfying the
following conditions:
\begin{enumerate}
\item[(i)] The set $\{\psi_0(\eps),\dots,\psi_k(\eps)\mid
    (\psi_0,\dots,\psi_k)\in\cb,\ \eps\in I \}$ is bounded in
    $\cd(\R^s)$;
\item[(ii)] $\sup\limits_{(\psi_i)\in\cb}\,\sup\limits_{i=0,\dots,s}
    |\lgl\xi^\bet,\psi_i(\eps)(\xi)\rgl|=O(\eps^q)$
    for all $\bet\in\N_0^s$ with $1\le|\bet|\le q$.
\end{enumerate}
\et
Note that $(\psi_0,\dots,\psi_k)\in \cc^\infty_b(I ,\ca_0(\R^s))\times
\big[\cc^\infty_b(I ,\ca_{00}(\R^s))\big]^k$ forms a $(k,q)$-class
$\{(\psi_0,\dots,\psi_k)\}$ (consisting of a single element) if and
only if each of $\psi_0,\dots,\psi_k$ has asymptotically vanishing
moments of order $q$.
The following results are established by combining techniques of the
respective proofs of Theorem 17 of \cite{JEL} and of Theorem
10.5.
\bt\lb{JT17A}
Let $R\in\ce(\Om)$. $R$ is moderate of type
$[\mathrm{A}_\mathrm{g}^\infty]$ if and only if the following
condition is satisfied:
$$\forall K\subset\subset\Om\ \forall\al\in\N_0^d\ 
   \forall k\in\N_0\ \exists N\in\N
   \ \mbox{such that for each $(k,N)$-class}\ \cb:
$$
$$\sup_{(\psi_i)\in\cb}\,\sup_{x\in K}|
    \pa^\al\rmd_1^kR_\eps(\psi_0(\eps),x)
    (\psi_1(\eps),\dots,\psi_k(\eps))|=O(\eps^{-N}).$$
Moreover, for given $K,\al,k,N$ the preceding condition is
satisfied for all $(k,N)$-classes $\cb$ if and only if it is satisfied
for all $(k,N)$-classes consisting of a single element
$(\psi_0,\dots,\psi_k)$. Therefore, the uniformity requirement
with respect to
$\cb$ can as well be omitted from the characterization of
moderateness given above.
\et
\bt\lb{JT1823A}
Let $R\in\ce(\Om)$. $R$ is negligible of type
$[\mathrm{A}_\mathrm{g}^\infty]$ if and only if the following
condition is satisfied:
$$\forall K\subset\subset\Om\ \forall\al\in\N_0^d\ 
   \forall k\in\N_0\ \forall n\in\N\ \exists q\in\N
   \ \mbox{such that for each $(k,q)$-class}\ \cb:
$$
$$\sup_{(\psi_i)\in\cb}\,\sup_{x\in K}|
\pa^\al\rmd_1^kR_\eps(\psi_0(\eps),x)
    (\psi_1(\eps),\dots,\psi_k(\eps))|=O(\eps^n).$$
Moreover, for given $K,\al,k,n,q$ the preceding condition is
satisfied for all $(k,q)$-classes $\cb$ if and only if it is satisfied
for all $(k,q)$-classes consisting of a single element
$(\psi_0,\dots,\psi_k)$. Therefore, the uniformity requirement
with respect to $\cb$ can as well be omitted from the characterization of
negligibility given above.
\et
The proofs of Theorems \rf{JT17A} and \rf{JT1823A} are deferred to
the end of this section.
\bc\lb{giki}
Let $R\in\ce(\Om)$. $R$ is moderate (resp.\ negligible)
with respect to type $[\mathrm{A}^\infty_{\mathrm{g}}]$ if and only if it is
moderate (resp.\ negligible)
with respect to type $[\mathrm{A}_\mathrm{l}^\infty]$.
\et
\pr
Sufficiency of the condition being obvious, let us show
necessity. Supposing $R$ to be moderate with respect to type
$[\mathrm{A}^\infty_{\mathrm{g}}]$, the differentials of $R$
satisfy the condition of Theorem \rf{JT17A}. For testing
$R$ on some $K\subset\subset\Om$ as to moderateness with
respect to type $[\mathrm{A}_\mathrm{l}^\infty]$, we have to 
consider test objects just of that type.
Now it is exactly the easy part of the very proof of \rf{JT17A}
which shows that this test gives a positive answer.
The same argument applies to negligibility.
\ep

From the preceding Corollary and Corollary \rf{ableitungengratis},
we see that all three types $[\mathrm{A}_\mathrm{g}^\infty]$,
$[\mathrm{A}_\mathrm{g}]$ and $[\mathrm{A}_\mathrm{l}^\infty]$
give rise to the same notions of moderateness resp.\ negligibility,
hence to the same Colombeau algebras. This fact also constitutes one of
the key ingredients for obtaining an intrinsic description of
the algebra $\cg^d(\Om)$ on manifolds: The property of a test object
living on the manifold to have asymptotically vanishing moments can
be formulated in intrinsic terms, indeed
(\cite{vi}, Definition 3.5); yet it would be virtually
unmanageable to deal with the latter property also for derivatives
of this test object, which, of course, are to be understood in this
general case as
appropriate Lie derivatives with respect to smooth vector fields.
Now Corollaries \rf{ableitungengratis} and \rf{giki}
allow to dispense with derivatives of test objects as regards the
asymptotic vanishing of the moments, provided all
$K\subset\subset\Om$ are taken into account (\cite{vi}, Corollary 4.5).

The next corollary might come as a bit of a surprise since we are
already used to type $[\mathrm{A}_\mathrm{l}]$ displaying rather bad
properties. Observe that it (necessarily, compare Example
7.7) only refers to moderateness. The case at hand
seems to be the only one where
a certain symmetry between $\ce_M$ and $\cn$ is
broken.
\bc\lb{gik}
Let $R\in\ce(\Om)$. $R$ is moderate
with respect to type $[\mathrm{A}_\mathrm{l}]$ if and only if it is
moderate
with respect to type $[\mathrm{A}_\mathrm{g}^\infty]$
(resp.\ $[\mathrm{A}_\mathrm{l}^\infty]$
resp.\ $[\mathrm{A}_\mathrm{g}]$).
\et
\pr       
Necessity of the condition being obvious this time,
let us show
sufficiency. Suppose $R$ to be moderate with respect to type
$[\mathrm{A}^\infty_{\mathrm{g}}]$ and let $K\subset\subset\Om$,
$\al\in\N_0^s$ be given. According to Theorem \rf{JT17A}, choose
$N\in\N$ such that for every $k=0,1,\dots,|\al|$, for every
$\bet\in\N_0^s$ with $0\le|\bet|\le|\al|$ and for every
$(k,N)$-class $\cb$,
$$\sup_{(\psi_i)\in\cb}\,\sup_{x\in K}|
    \pa^\al\rmd_1^kR_\eps(\psi_0(\eps),x)
    (\psi_1(\eps),\dots,\psi_k(\eps))|=O(\eps^{-N}).$$
For any test object $\phi$ of type
$[\mathrm{A}_\mathrm{l}]_{K,N(1+|\al|)}$,
it now follows
\beas
&&\!\!
\sup_K|\pa^\al(R_\eps(\phi(\eps,x),x))|\\
  &&\qquad=\sup_K\Big|\sum_{\bet,m}(\pa^\bet\rmd_1^mR_\eps)
    (\phi(\eps,x),x)
    (\pa^{\ga_1}\phi(\eps,x),\dots,\pa^{\ga_m}\phi(\eps,x))\Big|\\
  &&\qquad=
  \sup_K\Big|\sum_{\bet,m}(\pa^\bet\rmd_1^mR_\eps)
    (\phi(\eps,x),x)
    (\eps^N\pa^{\ga_1}\phi(\eps,x),\dots,\eps^N\pa^{\ga_m}\phi(\eps,x))
       \cdot\eps^{-mN}\Big|\\
  &&\qquad= O(\eps^{-N-|\al|N})
\eeas
since, for every $m\in\N_0$, the finite sequences
$(\phi(\eps,x),
\eps^N\pa^{\ga_1}\phi(\eps,x),\dots,\eps^N\pa^{\ga_m}\phi(\eps,x))$
(with $x$ ranging over $K$) form an $(m,N)$-class.
\ep

Now the proofs of {\bf (T4)} and {\bf (T5)},
that is, of the invariance of 
$\ce_M[\mathrm{A}_\mathrm{g}^\infty]$ and
$\cn[\mathrm{A}_\mathrm{g}^\infty]$ with respect to
differentiation, follow
from Theorems \rf{JT17A} and \rf{JT1823A}
in precisely
the same way as they have been achieved in section 7
for the builiding blocks of $\cg^d(\Om)$
by means of Theorems 7.12 and 7.13.
Digressing once more from the proof of
$\cg^2(\Om)$ being a diffeomorphism invariant Colombeau algebra,
let us deal with invariance under differentiation for the
remaining types of algebras:
Types $[\mathrm{A}_{\mathrm g}]$ and
$[\mathrm{A}_\mathrm{l}^\infty]$ as well as the case of
$\ce_M[\mathrm{A}_\mathrm{l}]$ are covered by
Corollaries \rf{ableitungengratis}, \rf{giki} and
\rf{gik}, respectively.
Moreover, it easy to check that {\bf (T4)} and {\bf (T5)}
are true for all types
$[\eps]$ and $[\mathrm{c}]$. An appropriate modification of
Theorem 17 of \cite{JEL} putting $\ca_N(\R^s)$
resp.\ $\ca_q(\R^s)$ in place of
$\ca_0(\R^s)$ and employing the techniques used in the proof of
$(\mathrm{A})\!\Leftrightarrow\!(\mathrm{C})$ of Theorem 10.5
establishes the respective results
to hold also for type $[\eps x,V]$.
Type $[\eps x,0]$ being covered by
section 7,
we are left only with $\cn[\mathrm{A}_\mathrm{l}]$ to be discussed.
However, lacking an
analog of
Theorem 7.13 resp.\ of Theorem \rf{JT1823A}
for type $[\mathrm{A}_\mathrm{l}]$ we are not in a position to
express negligibility
of $R$ with respect to this type in terms of differentials
of $R_\eps$. This tool, however, was the basis for deducing
invariance under differentiation. So for the time being, 
we find $\cn[\mathrm{A}_\mathrm{l}]$ to be
the only one among all 11+8 (to be precise, 8+6 pairwise different)
spaces $\ce_M[X]$ resp.\ $\cn[X]$ for
which the invariance with respect to differentiation
has to remain an open problem.

Finally, let us consider the question of diffeomorphism invariance.
As an inspection of the structure of the proof of
Theorem 7.14 of section 7
reveals, this theorem actually shows the
$\mu$-transform of test objects of types $[\mathrm{A}_\mathrm{l}]$,
$[\mathrm{A}_\mathrm{g}]$, $[\mathrm{A}_\mathrm{l}^\infty]$,
$[\mathrm{A}_\mathrm{g}^\infty]$
to be of the same type again,
respectively, thereby establishing {\bf (T6)}
in all four cases.
Moreover, we see that
on the basis of Corollaries \rf{ableitungengratis} and \rf{giki}
even a weaker version of the last statement of
Theorem 7.14, referring
only to types $[\mathrm{A}_\mathrm{l}]$ and
$[\mathrm{A}_\mathrm{g}]$, would suffice to obtain diffeomorphism
invariance for all four types $[\eps x,\mathrm{A}]$
(and, still, for $\ce_M[\eps x,0]$;
hence this would completely satisfy also the needs of
section 7 dealing with $\cg^d(\Om)$!):
Derivatives $\pa^\al_x\phi$ ($\al\neq0$) could be
dispensed with in Theorem 7.14 and its proof.

Recall that our proofs of {\bf (T7)} and {\bf (T8)}
in section 7 
(stating the invariance of moderateness resp.\ negligibility
under the action induced by a diffeomorphism)
were based on the equivalence
of conditions (C) and (Z) in Theorem 10.5
(resp.\ of $(\mathrm{C}'')$ and $(\mathrm{Z}'')$ in Corollary
10.7) which, in turn,
used the extension of paths $\phi(\eps,x)$ provided by
Proposition 10.4. Now each of the four
types $[\eps x,\mathrm{A}]$ is preserved by the extension process
$\phi\mapsto\ti\phi$.
Thus the respective analogs of
$(\mathrm{C})\!\Leftrightarrow\!(\mathrm{Z})$
can be shown to hold for all types $[\eps x,\mathrm{A}]$
by the methods employed in 10.5.

Now the proofs of {\bf (T7)} and {\bf (T8)}, respectively,
are literally the same for all four types $[\eps x,\mathrm{A}]$
as for the algebra $\cg^d(\Om)$
treated in section 7. The proof of {\bf (T8)} is even simpler
in the present case
since we do not have to bother with bridging the gap between
vanishing moments (as used in Definition 7.3) and
asymptotically vanishing moments (as occurring in Theorem 7.14)
which has been accomplished in section 7 by
means of the equivalence $(3^\circ)\!\Leftrightarrow\!(4^\infty)$
provided by Theorem 7.9.

Summarizing, we have established 
\bt\lb{thg2}
$\cg^2(\Om)$ is a
diffeomorphism invariant Colombau algebra which can be obtained by
using test objects of any of the types $[\mathrm{A}_\mathrm{g}^\infty]$,
$[\mathrm{A}_\mathrm{g}]$ or $[\mathrm{A}_\mathrm{l}^\infty]$.
\et
Test objects of type $[\mathrm{A}_\mathrm{l}]$, on the other hand,
give rise to a diffeomorphism invariant algebra
which does not preserve the product of smooth functions via $\io$
and for which it remains open
if it is a differential algebra at all.
Moreover, we have shown that each of the remaining six algebras
(apart from $\cg^d(\Om)=\cg[\eps x,0]$) satisfies {\bf (T1)}--{\bf
(T5)}, yet fails to be diffeomorphism invariant.

Also for the algebra $\cg^2(\Om)$ it is true that in characterizing
the negligibility of $R\in\ce_M(\Om)$ in terms of the differentials
of $R_\eps$, derivatives can be dispensed with,
those with respect to $\vphi$ as well as those with respect to $x$.
The numbering of the conditions in the following theorem
corresponds to that of Theorem \rf{thcond0}.
\bt\lb{thcond0g2}
Let $R\in\ce(\Om)$ be moderate with respect to type
$[\mathrm{A}_\mathrm{g}^\infty]$ (resp.\ with respect to types
$[\mathrm{A}_\mathrm{g}]$,
$[\mathrm{A}_\mathrm{l}^\infty]$).
Then $R$ is negligible with respect to any one of these
types if and only if one of the following (equivalent) conditions
is satisfied:\ms
$(0^\circ_\mathrm{A})$\quad
$\forall K\subset\subset\Om\  
   \forall n\in\N\ \exists q\in\N
   \ \mbox{such that for each $(0,q)$-class}\ \cb
$:
$$\sup_{(\psi_0)\in\cb}\,\sup_{x\in K}|
R_\eps(\psi_0(\eps),x)|
    =O(\eps^n).$$
$(1^\circ_\mathrm{A})$\quad
$\forall K\subset\subset\Om\ \forall\al\in\N_0^d\ 
   \forall n\in\N\ \exists q\in\N
   \ \mbox{such that for each $(0,q)$-class}\ \cb
$:
$$\sup_{(\psi_0)\in\cb}\,\sup_{x\in K}|
\pa^\al R_\eps(\psi_0(\eps),x)|
    =O(\eps^n).$$
$(2^\circ_\mathrm{A})$\quad
$\forall K\subset\subset\Om\ \forall\al\in\N_0^d\ 
   \forall k\in\N_0\ \forall n\in\N\ \exists q\in\N
   \ \mbox{such that for each $(k,q)$-class}\ \cb
$:
$$\sup_{(\psi_i)\in\cb}\,\sup_{x\in K}|
\pa^\al\rmd_1^kR_\eps(\psi_0(\eps),x)
    (\psi_1(\eps),\dots,\psi_k(\eps))|=O(\eps^n).$$
In each of the preceding conditions,
the uniformity requirement with respect to
$\cb$ can as well be omitted without changing the content of the
condition, regardless of the moderateness of $R$.
\et
\pr Due to Corollaries \rf{ableitungengratis} and \rf{giki}
it does not matter which of the three types is
being considered.
For $R\in\ce_M[\mathrm{A}_\mathrm{g}^\infty]$,
$(2^\circ_\mathrm{A})$ is equivalent to negligibility
with respect to $[\mathrm{A}_\mathrm{g}^\infty]$
by Theorem \rf{JT1823A}.
$(2^\circ_\mathrm{A})\!\Rightarrow\!(1^\circ_\mathrm{A})
\!\Rightarrow\!(0^\circ_\mathrm{A})$
is trivial; $(0^\circ_\mathrm{A})\!\Rightarrow\!(1^\circ_\mathrm{A})$ and
$(1^\circ_\mathrm{A})\!\Rightarrow\!(2^\circ_\mathrm{A})$
can be established by carefully
replacing bounded subsets of $\ca_0(\R^s)$ resp.\ of
$\ca_{00}(\R^s)$ by appropriately chosen
$(k,q)$-classes in the respective
proofs of Theorem \rf{thcond0} and
part $(1^\circ)\!\Rightarrow\!(2^\circ)$ of Theorem 
18 of \cite{JEL}.
As far as the proof of $(1^\circ_\mathrm{A})
\!\Rightarrow\!(2^\circ_\mathrm{A})$
(proceeding by induction with respect to $k$)
is concerned, the most delicate task in this respect consists in
choosing appropriate $(k\!+\!1,q)$- resp.\ $(k\!-\!1,q)$-classes
$\cb_{+1},\cb_{-1}$ to be used in connection with
$\pa^\al\rmd_1^{k+1}R_\eps$ resp.\ $\pa^\al\rmd_1^{k-1}R_\eps$
when $\pa^\al\rmd_1^kR_\eps$ is being evaluated on some
$(k,q)$-class $\cb$. To this end, define

$(k\!+\!1,q)$- resp.\ $(k\!-\!1,q)$-classes $\cb_{+1},\cb_{-1}$
by
\beas
\cb_{+1}&:=&\{(\psi_0+t\psi_k,\psi_1,\dots,\psi_k,\psi_k)\mid
   (\psi_0,\dots,\psi_k)\in\cb,\ 0\le t\le 1\},\\
\cb_{-1}&:=&\{(\psi_0+t\psi_k,\psi_1,\dots,\psi_{k-1})\mid
   (\psi_0,\dots,\psi_k)\in\cb,\ 0\le t\le 1\}
   \eeas
to provide the appropriate arguments for
$\pa^\al\rmd_1^{k+1}R_\eps$ resp.\
$\pa^\al\rmd_1^{k-1}R_\eps$.
On the basis of this choice of
$\cb_{+1},\cb_{-1}$ the proof of Theorem 18 of \cite{JEL} can be
upgraded by introducing $\eps$ as additional parameter throughout
as to establish
$(1^\circ_\mathrm{A})\!\Rightarrow\!(2^\circ_\mathrm{A})$
of the theorem.

Note that in the proof of
$(0^\circ_\mathrm{A})\!\Rightarrow\!(1^\circ_\mathrm{A})$
being obtained from the proof of
Theorem \rf{thcond0} by introducing the parameter $\eps$,
Theorem {\bf (T4)} which,
in turn, is based on Theorem \rf{JT17A} has to be
invoked
to guarantee the moderateness of $\pa_iR$.

The last statement of the theorem
follows from Theorem \rf{JT1823A} since the corresponding
statement thereof
contains, among others, $\al$ and $k$ as free variables.
\ep

A glance at the preceding proof shows virtually
all the substantial results of this article to be involved.
Note that a $(0,q)$-class consisting of a single element
$\phi$ is nothing else than (a singleton containing)
a test object of type
$[\eps,\mathrm{A}]_q$.
The equivalence of $R\in\cn[\mathrm{A}_\mathrm{g}^\infty]$
and condition $(0^\circ_\mathrm{A})$ without the uniformity
clause (provided $R\in\ce_M[\mathrm{A}_\mathrm{g}^\infty]$)
now shows that for a function $R\in\ce(\Om)$
which is moderate with respect to type
$[\eps x,\mathrm{A}_\mathrm{g}^\infty]$, it amounts to the same
to be negligible with respect to either 
type $[\eps x,\mathrm{A}_\mathrm{g}^\infty]$
or $[\eps,\mathrm{A}]$. We will make use of this fact below.
For the following theorem, recall that $\cg^e_0$ denotes the smooth
part of $\cg^e$ (cf.\ the discussion following Definition
\rf{colalg}).

\bt\lb{injcanon}
Of the canonical maps
$\cg^d(\Om)\to\cg^2(\Om)\to\cg^1(\Om)\to\cg^e_0(\Om)$
the first and the second one are injective
whereas the third one is not. 
The four corresponding spaces of representatives (i.e., of moderate
functions) are pairwise different.
\et
\pr
The injectivity of the map $\cg^d(\Om)\to\cg^2(\Om)$
is equivalent to
$\cn[\eps x,\mathrm{V}]\cap\ce_M[\eps x,0] 
   =\cn[\eps x,\mathrm{A}_\mathrm{g}^\infty]\cap\ce_M[\eps x,0]$.
This, however, is accomplished by the extension
$(3^\circ)\!\Leftrightarrow\!(5^\infty)$ 
of Theorem 7.9
derived from the first diagram in section \rf{spec}.
The injectivity of $\cg^2(\Om)\to\cg^1(\Om)$, on the other hand,
is equivalent to
$\ce_M[\eps x,\mathrm{A}_\mathrm{g}^\infty]
   \cap\cn[\eps,\mathrm{A}] 
   =\cn[\eps x,\mathrm{A}_\mathrm{g}^\infty]$
which has been deduced previously from Theorem~\rf{thcond0g2}.
Finally, to establish the non-injectivity of $\cg^1(\Om)\to\cg^e_0(\Om)$,
we use
the fact that the function $P$ introduced in section \rf{ex}
can be shown not to be negligible with respect to type
$[\eps,\mathrm{A}]$, by techniques similar to those
employed in section \rf{ex}.
The difference of the respective spaces of moderate functions
should be clear from the following examples.
\ep

\bexs\lb{R2345}
Let $\Om:=\R$.
\begin{itemize}
\item[(i)] Both $R_2(\vphi,x):=\exp(\lgl\vphi|\vphi\rgl^2
   \lgl\xi,\vphi(\xi)\rgl)$ and
   $R_3(\vphi,x):=\exp(\vphi(0)^2
   \lgl\xi,\vphi\rgl)$ are moderate of type
   $[\eps x,\mathrm{A}^\infty_{\mathrm{g}}]$, yet not of type
   $[\eps x,0]$.
\item[(ii)] The counterexample
   $R_0(\vphi,x):=\exp(i\exp(\lgl\vphi|\vphi\rgl))$
   (see \cite{JEL}) which has already been mentioned in
   connection with the failure of algebras of type $[\eps]$
   to be diffeomorphism invariant (see \rf{R01}(i))
   has the property
   of being moderate of type $[\eps,\mathrm{A}]$
   yet not of type $[\eps x,\mathrm{A}^\infty_{\mathrm{g}}]$.
   The same holds true for $R_5$ to be defined below.
\item[(iii)] $R_4(\vphi,x):=\lgl\xi,\vphi(\xi)\rgl
   \cdot\exp(\vphi(0))$ is moderate of type $[\mathrm{c},\mathrm{V}]$,
   yet not of type $[\eps,\mathrm{A}]$. This also holds true for
   $R_1$ introduced as Example \rf{R01}\,(ii).
\item[(iv)] $R_5(\vphi,x):=\exp(-\lgl\vphi|\vphi\rgl)\cdot
   \exp(i\exp(2\lgl\vphi|\vphi\rgl))$
   is of particular interest: It is not moderate with respect to
   any of the types $[\eps x]$, however, it is even negligible
   with respect to all types $[\eps]$ and $[\mathrm{c}]$.
\end{itemize}
Again the proofs of the preceding claims are left to the reader.
\et
It is a remarkable fact that answering the
apparently harmless question of injectivity of
the canonical maps in the last analysis involves
quite a number of hard theorems: the extension 
$(3^\circ)\!\Leftrightarrow\!(5^\infty)$
of Theorem 7.9
derived in section \rf{spec}; Theorem \rf{thcond0g2} which, in turn,
is based on part
$(1^\circ)\!\Leftrightarrow\!(2^\circ)$
Theorem 18 of \cite{JEL} and on Theorems \rf{thcond0},
\rf{JT17A} and \rf{JT1823A}; finally, also
the counterexample $P$ of section \rf{ex} is among the ingredients
of the argument. It remains to prove Theorems
\rf{JT17A} and \rf{JT1823A}.

{\bf Proof of Theorem \rf{JT17A}. }
To show sufficiency of the condition, suppose that the
differentials of $R_\eps$ (where
$R\in\ce(\Om)$) satisfy the property
specified in the theorem. Consider a test object
$\phi(\eps,x)$ of type $[\mathrm{A}_\mathrm{g}^\infty]_N$ and set
$\Phi(\eps,x):=(\phi(\eps,x),x)$. Expanding
$\pa^\al(R_\eps\circ\Phi)$ according to the chain rule shows that
$R$ is moderate of type $[\mathrm{A}_\mathrm{g}^\infty]$:
It suffices to observe that the family of all finite sequences
$$(\phi(\eps,y), \pa_y^{\bet_1}\phi(\eps,y),
   \dots,\pa_y^{\bet_l}\phi(\eps,y))$$
forms an $(l,N)$-class if $\eps$ is considered as variable and $y$
as a parameter taking values in some compact subset of $\Om$.

Conversely, for a function $R\in\ce(\Om)$ which is moderate with
respect to type $[\mathrm{A}_\mathrm{g}^\infty]$, we will show
that the assumption of $R$ to violate the condition in the theorem
leads to a contradiction. Thus suppose that there exist
$K\subset\subset\Om$, $\al\in\N_0^s$, $k\in\N_0^s$ such that for
all $N\in\N$ there exists a $(k,N)$-class $\cb$ such that
\bea\lb{1*}
\sup_{K,\cb}|\pa^\al\rmd_1^k R_\eps(\psi_0(\eps),x)
   (\psi_1(\eps),\dots,\psi_k(\eps))|
\eea
is {\it not} of order $\eps^{-N}$. By moderateness of $R$, there
exists $N\in\N$ such that
\bea\lb{2*}
\sup_K|\pa^{\al'}(R_\eps(\phi(\eps,x),x))|=O(\eps^{-N})
\eea
for all test objects $\phi$ of type
$[\mathrm{A}_\mathrm{g}^\infty]_N$, where
$\al':=\al+pe_s$, $p:=\sum\limits_{i=1}^{k}(|\al|+k^2+i)$.
\linebreak
Due to our hypothesis, there exists a $(k,N)$-class $\cb$
such that (\rf{1*}) is not of order $\eps^{-N}$.
Having fixed $K,\al,k,N,\cb$, we inductively
define sequences $x^{(j)}\in K$,
$(\psi_0^{(j)},\dots,\psi_k^{(j)})\in\cb$,
$0<\eps_j<\frac{1}{j}$ (with $\eps_{j+1}<\eps_j$) ($j=1,2,\dots$)
such that the following inequalities hold for $j=1,2,\dots$:
\bea\lb{3*}
|\pa^\al\rmd_1^k R_{\eps_j}(\psi_0^{(j)}(\eps_j),x^{(j)})
   (\psi_1^{(j)}(\eps_j),\dots,\psi_k^{(j)}(\eps_j))|
   \ge j\cdot\eps_j^{-N}
\eea
(the technical details are similar to those in the proof of part
(C)$\Rightarrow$(A) of Theorem 10.5). Let $(\lambda_j)_{j\in\N}$
be a partition of unity on $I $ as in Lemma 10.1;
for $(t_1,\dots,t_k)\in\{0,1,\dots,k\}^k$, define
$$
\phi_{t_1,\dots,t_k}(\eps,x):=
   \sum_{j=1}^{\infty} \lambda_j(\eps)\cdot
   \bigg[\psi_0^{(j)}(\eps)+
   \sum_{i=1}^{k}t_i\frac{(x_s-x_s^{(j)})^{|\al|+k^2+i}}
   {(|\al|+k^2+i)!}
   \cdot\psi_i^{(j)}(\eps)\bigg].
$$
Since $\sum\limits_{i=1}^{k}t_i\frac{(x_s-x_s^{(j)})^{|\al|+k^2+i}}
   {(|\al|+k^2+i)!}$ is a polynomial in $x$ and all
$(\psi_0^{(j)},\dots,\psi_k^{(j)})$ are members of one
particular $(k,N)$-class $\cb$, $\phi_{t_1,\dots,t_k}$
is a member of $\cc^\infty_b(I ,\ca_0(\R^s))$ and, in addition,
is of type $[\mathrm{A}_\mathrm{g}^\infty]_N$.
From (\rf{2*}) we conclude that
\bea\lb{4*}
\sup_K|\pa^{\al'}(R_\eps(\phi_{t_1,\dots,t_k}(\eps,x),x))|
   =O(\eps^{-N}).
\eea
Now we follow the combinatorial reasoning
of the proof of Theorem 17 of \cite{JEL} to derive the desired
contradiction: Choosing numbers $c_0,\dots,c_k$ satisfying
the set of equations $\sum\limits_{i=0}^{k}c_i\cdot i^m=\de_{1m}$
($m=0,1,\dots,k$), let us form
\bea\lb{5*}
\sum_{t_1=0}^{k}\ldots\sum_{t_k=0}^{k}c_{t_1}\dots c_{t_k}
   \pa^{\al'}(R_\eps(\phi_{t_1,\dots,t_k}(\eps,x),x)).
\eea
By (\rf{4*}), this expression is of order $\eps^{-N}$,
uniformly for $x\in K$. On the other hand, evaluating (\rf{5*})
at $\eps:=\eps_j$, $x:=x^{(j)}$ according to the chain rule results
in a positive integer multiple of
$\pa^\al\rmd_1^k R_{\eps_j}(\psi_0^{(j)}(\eps_j),x^{(j)})
   (\psi_1^{(j)}(\eps_j),\dots,\psi_k^{(j)}(\eps_j))$,
due to the delicate combinatorial argument of the proof of Theorem
17 of \cite{JEL}. (\rf{5*}) being of order $\eps^{-N}$, we conclude
that
\beas
|\pa^\al\rmd_1^k R_{\eps_j}(\psi_0^{(j)}(\eps_j),x^{(j)})
   (\psi_1^{(j)}(\eps_j),\dots,\psi_k^{(j)}(\eps_j))|
   \le C\eps_j^{-N}\qquad (j\ge j_0)
   \eeas
for some positive constant $C>0$ and some $j_0\in\N$. This,
however, contradicts our choice of $x^{(j)}$,
$\psi_i^{(j)}$. So the condition in the theorem, in fact, is
necessary for $R$ being moderate.
(In a trivial way, the preceding reasoning also
applies in the case $k=0$ if all sums $\sum\limits_{i=1}^{k}$
are set equal to $0$.)

To prove the last statement of the theorem,
let $K,\al,k,N$ be given. Suppose, again by way of
contradiction, the condition on the differentials of
$R_\eps$ given in the theorem
to be satisfied for $(k,N)$-classes consisting of a single
element, yet to be violated for arbitrary $(k,N)$-classes,
either with resepect to the particular $K,\al,k,N$ at hand.
Similarly to the reasoning of the main part of the proof, deduce
from these hypotheses the existence of
a $(k,N)$-class $\cb$ and of sequences
$0<\eps_{j+1}<\eps_j<\frac{1}{j}$, $x^{(j)}\in K$,
$(\psi_0^{(j)},\dots,\psi_k^{(j)})\in\cb$ ($j=1,2,\dots$)
satisfying the inequalities (\rf{3*}) for all $j\in\N$.
Now define
$$\psi_i(\eps):=\sum_{j=1}^{\infty}\lambda_j(\eps)\psi_i^{(j)}(\eps)
    \qquad\qquad(i=0,1,\dots,k)$$
where $(\lambda_j)_{j\in\N}$
is a partition of unity on $I $ as in Lemma 10.1.
Due to the properties of the $\psi_i^{(j)}$,
$\{(\psi_0,\dots,\psi_k)\}$ is a $(k,N)$-class.
By assumption,
$$\sup_K|\pa^\al\rmd_1^k R_\eps(\psi_0(\eps),x)
   (\psi_1(\eps),\dots,\psi_k(\eps))|=O(\eps^{-N}).$$
Taking into account that $\psi_i(\eps_j)=\psi_i^{(j)}(\eps_j)$, this
contradicts our choice of $x^{(j)}$, $\psi_i^{(j)}$, thereby
completing the proof.
\ep

{\bf Proof of Theorem \rf{JT1823A}. }
Just copy the proof of Theorem \rf{JT17A}, add ``for all $n\in\N$''
at the appropriate places and change ``$\eps^{-N}$''
to ``$\eps^n$''. At the remaining occurrences of $N$, replace it by
$q$.
\ep

\section{
Concluding remarks}\lb{conclu}
As has been pointed out already in Part I, Theorem
7.14 has a place at the very core of
diffeomorphism invariance of a Colombeau algebra.
The problem with algebras of any type $[\mathrm{c},Y]$, of course,
is that classes consisting of test objects as simple as
$\vphi\in\ca_0(\R^s)$ are not
invariant under the action of a diffeomorphism since the latter
introduces dependence on $\eps$ and $x$. Types $[\eps x,Y]$
of course are an efficient remedy against that problem as they
incorporate a very general $(\eps,x)$-dependence into test objects.
Yet there is an intermediate way: Starting with the class of
``constant'' test
objects $\ti\vphi\in\ca_0(\R^s)$ resp.\ $\in\ca_q(\R^s)$, we consider
the minimal class containing these which is invariant with respect
to diffeomorphisms. Due to the functorial property of $\bar\mu_\eps$,
this class precisely consists of all images of constant test objects (in the
sense just described) under the mappings
$\ti\vphi\mapsto((\eps,x)\mapsto\pro_1\bar\mu_\eps(\ti\vphi,\mu^{-1
}x))$
where $\mu$ ranges over all diffeomorphisms onto the open set under
consideration. As the following example shows, the class of test
objects obtained in this way (starting with all
$\ti\vphi\in\ca_0(\R^s)$) is in fact different in general from
$\cc^\infty_b(I \times\Om,\ca_0(\R^s))$.
\bex
Let $\Om:=\R$; choose $\psi\in\ca_0(\R)$ satisfying $\psi(0)\ne0$,
$\supp\psi\subseteq[-1,+1]$. Setting $\phi(\eps,x)(\xi):=
\psi(\xi+\sin x)$ in fact defines an element of
$\cc^\infty_b(I \times\R,\ca_0(\R))$. Now assume that there
exists $\ti\vphi\in\ca_0(\R)$ and a diffeomorphism $\mu:\ti\Om\to\R$
(where $\ti\Om\subseteq\R$ is open) such that
$(\phi(\eps,x),x)=\bar\mu_\eps(\ti\vphi,\mu^{-1}x)$ for all $x\in\R$.
Setting $\xi:=0$, $x_1:=0$, $x_2:=\frac{\pi}{2}$, respectively,
we obtain $\psi(0)=\ti\vphi(0)\cdot|(\mu^{-1})'(\mu(0))|$ resp.\
$\psi(1)=\ti\vphi(0)\cdot|(\mu^{-1})'(\mu(\frac{\pi}{2}))|$.
The first of these relations entails $\ti\vphi(0)\ne0$
while the second one implies $\ti\vphi(0)=0$, so we arrive at a
contradiction.
\et
In some situations, it may not even be desirable to require
invariance of a Colombeau algebra with respect to all
diffeomorphisms; for example, invariance only with respect to
members of the Poincar\'e group might be of interest in
applications in special relativity. This approach also could be
combined with restricting the class of test objects
to images of constant test objects under the particular group
of transformations at hand. This opens the way to new
classes of Colombeau algebras possessing weaker invariance properties
than $\cg^d(\Om)$ does.
However, these new objects still can be constructed on the basis of the scheme
outlined in section 3 which, in our view, constitutes an
appropriate
general framework for the treatment of (full) Colombeau algebras.

{\bf Acknowledgements.} 
The work on this series of papers was initiated during a visit of the authors
at the department of mathematics of the university of Novi Sad in July 1998. 
We would like to thank the faculty and staff
there, in particular Stevan Pilipovi\'c and
his group for many helpful 
discussions and for their warm hospitality. Also, we are indebted to Andreas 
Kriegl for sharing his expertise on infinite dimensional calculus.

{\small
{\it Electronic Mail:} {\tt michael@mat.univie.ac.at}
}

\end{document}